\documentclass[12pt,oneside,reqno]{amsart}
\usepackage{}
\usepackage{amssymb}
\usepackage[colorlinks=true, linkcolor=blue, citecolor=red]{hyperref}
\usepackage{bbm}
\usepackage{cases}
\usepackage{cite}
\usepackage{amsmath}
\usepackage{graphicx}
\usepackage{mathrsfs}
\usepackage{stmaryrd}
\usepackage{color}
\usepackage{soul}
\usepackage[dvipsnames]{xcolor}
\usepackage{amsfonts}
\usepackage{enumerate,amsmath,amssymb,amsthm}

\numberwithin{equation}{section}

\newcommand{\be}{\begin{eqnarray}}
\newcommand{\mE}{\end{eqnarray}}
\newcommand{\ce}{\begin{eqnarray*}}
\newcommand{\de}{\end{eqnarray*}}
\newtheorem{theorem}{Theorem}[section]
\newtheorem{lemma}[theorem]{Lemma}
\newtheorem{remark}[theorem]{Remark}
\newtheorem{definition}[theorem]{Definition}
\newtheorem{proposition}[theorem]{Proposition}
\newtheorem{example}[theorem]{Example}
\newtheorem{corollary}[theorem]{Corollary}

\def\eps{\varepsilon}

\def\t{\tau}

\def\a{\alpha}

\def\b{\beta}
\def\p{\partial}

\def\[{{\Big[}}
\def\]{{\Big]}}
\def\<{{\langle}}
\def\>{{\rangle}}
\def\({{\Big(}}
\def\){{\Big)}}

\def\bx{{\mathbf{x}}}
\def\tr{{\rm tr}}

\def\dif{{\mathord{{\rm d}}}}

\def\no{\nonumber}
\def\={&\!\!=\!\!&}
\def\bt{\begin{theorem}}
\def\et{\end{theorem}}
\def\bl{\begin{lemma}}
\def\el{\end{lemma}}
\def\br{\begin{remark}}
\def\er{\end{remark}}

\def\bd{\begin{definition}}
\def\ed{\end{definition}}
\def\bp{\begin{proposition}}
\def\ep{\end{proposition}}
\def\bc{\begin{corollary}}
\def\ec{\end{corollary}}
\def\bx{\begin{example}}
\def\ex{\end{example}}

\def\cI{{\mathcal I}}

\def\cL{{\mathcal L}}

\def\cT{{\mathcal T}}

\def\mE{{\mathbb E}}

\def\mI{{\mathbb I}}

\def\mN{{\mathbb N}}

\def\mR{{\mathbb R}}

\def\sH{{\mathcal H}}

\def\sL{{\mathcal L}}

\def\sR{{\mathcal R}}
\def\sS{{\mathcal S}}

\def\sZ{{\mathcal Z}}

\def\geq{\geqslant}
\def\leq{\leqslant}

\def\t{\mathord{{\rm tr}}}

\allowdisplaybreaks

\begin{document}

\title{Infinite-horizon time-inhomogeneous Kolmogorov equation and  non-autonomous multi-scale stochastic systems}

\date{}

\author{Ling Wang, Pengcheng Xia, Longjie Xie and Li Yang}

\address{Ling Wang:
	School of Mathematics and Statistics, Jiangsu Normal University,
	Xuzhou, Jiangsu 221000, P.R. China\\
	Email: lwangmath@jsnu.edu.cn
}
\address{Pengcheng Xia:	School of Mathematics and Big Data, Anhui University of Science and Technology, Huainan,
 Anhui 232001, P.R.China\\
 	Email: pcxia@whu.edu.cn
}

\address{Longjie Xie:
	School of Mathematics and Statistics, Jiangsu Normal University,
	Xuzhou, Jiangsu 221000, P.R. China\\
	Email: longjiexie@jsnu.edu.cn
}

\address{Li Yang:
	School of Mathematical Sciences, Tiangong University,
	Tianjin 300387, P.R.China\\
	Email: lyang@tiangong.edu.cn
}
\thanks{
	This work is supported by the National Key R$\&$D program of China (No. 2023YFA1010103) and NNSF of China (No.  12471140, 12301179,  12401182).
}

\begin{abstract}

We  propose a unified approach,
based on the infinite-horizon time-inhomogeneous Kolmogorov equation with a parameter, to study the asymptotic behavior of  non-autonomous multi-scale stochastic systems with irregular coefficients, where both the fast and the slow equations depend on the highly oscillating time component.  We establish  both the weak and strong convergence of double  averaging principle as well as the functional central limit theorem with new homogenized coefficients given in terms of the evolution system of invariant measures of the time-inhomogeneous fast frozen equation.  Moreover,  rates of convergence are obtained as byproducts of our argument, which do not depend on the   regularities of the coefficients with respect to the time component and the fast variable.
	
	\bigskip

	\noindent {{\bf AMS 2010 Mathematics Subject Classification:} 60H10, 60F17, 35B30.}
	
	\vspace{2mm}
	\noindent{{\bf Keywords and Phrases:} Time-inhomogeneous Kolmogorov equation; non-autonomous multi-scale system; averaging principle; central limit theorem.}
\end{abstract}

\maketitle

\section{Introduction}
Our goal in this paper is to establish a general result for the asymptotic behavior of the non-autonomous multi-scale stochastic system. Before introducing the non-autonomous model, let us first recall the following autonomous stochastic  differential equation (SDE for short)  in $\mR^{d_1+d_2}$:
\begin{equation} \label{sde-au}
\left\{ \begin{aligned}
&\dif X^{\eps}_t =\frac{1}{\eps}b(X^{\eps}_t, Y^{\eps}_t)\dif t+\frac{1}{\sqrt{\eps}}\sigma(X_t^\eps,Y_t^\eps)\dif W^{1}_t,\qquad\quad X_0^\eps=x\in\mR^{d_1},\\
&\dif Y^{\eps}_t =F(X^{\eps}_t,Y^{\eps}_t)\dif t+G(X_t^\eps,Y_t^\eps)\dif W^{2}_t,\qquad\qquad\quad Y_0^\eps=y\in\mR^{d_2},
\end{aligned} \right.
\end{equation}
where $0<\eps\ll 1$ is a small parameter which represents the separation of time scales between the slow process $Y_t^\eps$   and the fast motion $X_t^\eps$.
Such multi-scale models have a wide range of applications including climate-weather interactions  (see e.g. \cite{K,MTV}),  macro-molecules (see e.g. \cite{BKRP}), geophysical fluid flows (see e.g. \cite{GD}), stochastic volatility in finance (see e.g. \cite{FFK}), etc. However, it is often too difficult to analyze or simulate the underlying system (\ref{sde-au}) directly due to the widely separated time scales and the cross interactions between the slow and fast modes.  Then the celebrated theory of the averaging principle says that a good approximation of the slow component can be obtained by averaging   with respect to the fast variable. More precisely, under suitable  assumptions on the coefficients,
the slow motion $Y_t^\eps$ will converge weakly (i.e., in distribution)  as $\eps\to 0$ to  $\bar Y_t$ which satisfies  the following averaged equation in $\mR^{d_2}$:
\begin{align}\label{sde-bau}
\dif \bar Y_t=\bar F( \bar Y_t)\dif t+ \bar G(\bar Y_t)\dif W^2_t,\quad \bar Y_0=y\in\mR^{d_2},
\end{align}
where the new averaged coefficients are defined by
\begin{align}\label{Fy}
\bar F(y):=\int_{\mR^{d_1}}F(x,y)\mu^y(\dif x),\quad\bar{G}(y):=\sqrt{\int_{\mR^{d_1}}G(x,y)G(x,y)^{*}\mu^y(\dif x)},
\end{align}
and for each $y\in\mR^{d_2}$, $\mu^y(\dif x)$ is the unique invariant measure for $X_t^y$ which satisfies  the frozen equation
\begin{align}\label{sde-fau}
\dif X_t^y=b(X_t^y,y)\dif t+\sigma(X_t^y,y)\dif W_t^1,\quad X_0^y=x\in\mR^{d_1}.
\end{align}
Furthermore, when the diffusion coefficient $G$ in the slow process does not depend on the fast motion, i.e.,
\begin{align}\label{gg}
 G(x,y)\equiv G(y),
\end{align}
then the strong convergence (i.e., in $L^2(\Omega)$-sense) also holds (a counterexample exists showing that the strong convergence does not hold if $G$ depends on the $x$-variable, see e.g. \cite{LD}).
The reduced equation (\ref{sde-bau}) then captures the essential dynamics of the system (\ref{sde-au}), which does not depend on the fast variable anymore and thus is much simpler.
Such results have been intensively studied for various multi-scale stochastic systems in the past decades, see  \cite{HL,Kh,RX} and the references therein.

The strong convergence in the averaging principle can be regarded as a functional law of large numbers.  In this case,
fluctuations of $Y_t^\eps$ around its
averaged motion $\bar Y_t$ have also been studied, see e.g. \cite{RX1,RXY,WR}, and it was shown that the normalized difference
$$
Z_t^\eps:=\frac{Y_t^\eps-\bar Y_t}{\sqrt{\eps}}
$$
will converge weakly   to a Gaussian process   which satisfies
\begin{align}\label{bzee}
\dif \bar Z_{t}=\nabla_y\bar{ F}(\bar Y_t)\cdot\bar Z_{t}\dif t+\nabla _y  G(\bar Y_t)\cdot\bar Z_{t}\dif W_t^2+\Gamma(\bar Y_t)\dif \tilde W_t,\quad \bar Z_0=0,\end{align}
where $\bar{ F}(y)$ is defined in (\ref{Fy})  ,  $\tilde W_t$ is another Brownian motion independent of $W_t^2,$
and the new homogenized diffusion coefficient is given by
\begin{align*}
\frac{1}{2}\Gamma\Gamma^*(y):=\int_{\mR^{d_1}}\!\big(F(x,y)-\bar F(y)\big)\Phi^*(x,y)\mu^y(\dif x),
\end{align*}
with $\Phi(x,y)$ being the unique solution of the following  Poisson equation in $\mR^{d_1}$:
\begin{align}\label{poi}
\sL_0(x,y)\Phi(x,y)=-\big(F(x,y)-\bar F(y)\big),
\end{align}
where $y\in\mR^{d_2}$ is a parameter, and the operator $\sL_0(x,y)$ is the generator of the ergodic homogeneous  process $X_t^y$ satisfying SDE (\ref{sde-fau}). Such a result is an analogue of the classical  central limit theorem and is
closely related to the homogenization of solutions of the second-order elliptic or parabolic equations, see, e.g. \cite{HP,Par}.


It is worth  pointing  out that the Poisson equation (\ref{poi})  not only appears in the limit equation (\ref{bzee}), but also serves as  one of the main tools in establishing the  averaging principle and the central limit theorem.  This equation is also referred to as the cell problem in the context of periodic
homogenization. Besides, it proves to be a powerful tool in diffusion approximation,  moderate and large deviations, numerical approximation of time-averaged estimators and the invariant measure of homogeneous SDEs or SPDEs,  and other limit theorems in probability theory, see e.g. \cite{NX,MST,P-V2,RX,SX} and the references therein.

However, there are  fewer results concerning the asymptotic behavior of the non-autonomous multi-scale systems. Namely, consider the following SDE  in $\mR^{d_1+d_2}$:
\begin{equation} \label{sde0}
\left\{ \begin{aligned}
&\dif X^{\eps}_t =\frac{1}{\eps}b(t/\eps, X^{\eps}_t, Y^{\eps}_t)\dif t+\frac{1}{\sqrt{\eps}}\sigma(t/\eps,X_t^\eps,Y_t^\eps)\dif W^{1}_t,\quad\, X^{\eps}_0=x\in\mR^{d_1},\\
&\dif Y^{\eps}_t =F(t/\eps,X^{\eps}_t,Y^{\eps}_t)\dif t+G(t/\eps,X^{\eps}_t,Y_t^\eps)\dif W^{2}_t,\;\quad\qquad\! Y^{\eps}_0=y\in\mR^{d_2}.
\end{aligned} \right.
\end{equation}
The novelty of  system (\ref{sde0}) lies in the fact that both the fast and the slow equations
include  the highly oscillating time component $t/\eps$.
Such systems,  with an external time-dependent perturbation, arise naturally  in many areas of biology and physics. Typical examples are climate-weather interactions involving diurnal cycle and seasonal cycle (see e.g. \cite{BM,CSG}),
and  neural networks with  time-dependent inputs (synaptic activities) (see e.g. \cite{GW}).
Yet the analysis of the non-autonomous system (\ref{sde0}) does not follow in a straightforward way from those results of autonomous systems available in the   literature.

When only the slow process involves the highly oscillating time component while the fast motion does not depend on the time variable, i.e., $b(t,x,y)\equiv b(x,y)$ and $\sigma(t,x,y)\equiv\sigma(x,y)$ in (\ref{sde0}), the averaging principle was studied in \cite{CL} (see also \cite{YSW} for similar results),
where the averaged drift  is shown to be related to the KBM-type vector field (KBM stands for Krylov, Bogolyubov and Mitropolsky, see \cite[Definition 4.2.4]{SVM}).  The situation is even more different when the highly oscillating time component appears in the fast motion. In this case, the frozen equation corresponding to the non-autonomous system (\ref{sde0}) should be chosen as the following time-inhomogeneous SDE: for $t\geq s>0$,
\begin{align}\label{froz}
\dif X_{s,t}^y=b(t,X_{s,t}^y,y)\dif t+\sigma(t,X_{s,t}^y,y)\dif W_t^1,\quad X_s^y\in\mR^{d_1},
\end{align}
where $y\in\mR^{d_2}$ is a parameter. Unlike the previous frozen equation (\ref{sde-fau}), we cannot hope to have a single invariant measure for the inhomogeneous Markov process $X_{s,t}^y(x)$ since the coefficients of (\ref{froz}) depend on the time variable.
Instead, we need to look for an evolution system of invariant measures that are families of  probability measures
$\{\mu_t^y\}_{t\geq 0}$ such that
\begin{align}\label{ergo}\int_{\mR^{d_1}}P^y_{s,t}\varphi(x)\mu_s^y(\dif x)=\int_{\mR^{d_1}}\varphi(x)\mu_t^y(\dif x),\quad\forall s\leq t,\varphi\in C_b(\mR^{d_1}),\end{align}
where $P^y_{s,t}$ is the two-parameter  semigroup associated with SDE
(\ref{froz}), i.e., for every $s\leq t$ and $x\in\mR^{d_1}$,
$$
P_{s,t}^y\varphi(x):=\mE\varphi\big(X^y_{s,t}(x)\big).
$$
Under suitable dissipative assumptions (in particular, if the drift
is strongly confining),  the existence of the evolution system of invariant measures $\{\mu_t^y\}_{t\geq 0}$ satisfying (\ref{ergo}) for  the frozen equation  (\ref{froz}) has been studied in  \cite{DR1,FZ,LLZ}, see also \cite{JQ,KLL} and the references therein.
When the slow process  does not depend on the time variable, i.e., $F(t,x,y)\equiv F(x,y)$ and $G(t,x,y)\equiv 0$ in (\ref{sde0}),  and the fast motion involves the highly oscillating time component and is $\tau$-periodic in time, the corresponding averaging principle was derived in \cite{Wg} (see also \cite{GW,Uda} for further extensions), namely, $Y_t^\eps$ converges strongly as $\eps\to0$ to $\bar Y_t$ which satisfies
\begin{align*}
\dif \bar Y_t=\bar{\bar F}(\bar Y_t)\dif t,\quad\bar Y_0=y\in\mR^{d_2},
\end{align*}
where the averaged drift is defined by
\begin{align*}
\bar{\bar F} (y):=\frac{1}{\tau}\int_0^{\tau}\!\!\!\int_{\mR^{d_1}}F(x,y)\mu_t^y(\dif x)\dif t.
\end{align*}
The case of a uniformly almost periodic fast system (while the slow equation does not involve the fast oscillating time component) was  considered for stochastic reaction-diffusion equations in \cite{Ce}. We remark that the main method used in the above works is the classical Khasminskii's time discretization  argument, and no rates of convergence are obtained. Moreover, to the best of our knowledge, there is still no result concerning  the normal deviations for the non-autonomous multi-scale systems.

We introduce a new approach,
based on the infinite-horizon time-inhomogeneous Kolmogorov equation, to study  the  averaging principle as well as  the  normal deviations of the non-autonomous multi-scale system (\ref{sde0}) with irregular coefficients.
More precisely, we  first establish a double averaging type principle for  (\ref{sde0}), including  both the weak and strong convergence, see {\bf Theorem \ref{main1}} and {\bf Theorem \ref{main11}} below.
When the strong convergence  holds,  we  further study the  fluctuations of $Y_t^\eps$ around its
average $\bar Y_t$. Define the normalized deviation process as
\begin{align}\label{zep}
Z_t^\eps:=\frac{Y_t^\eps-\bar Y_t}{\eps^{(1/2)\wedge\vartheta}},
\end{align}
where  $\vartheta$ is given in assumption ${\bf ({\tilde H_2})}$ below. It turns out that, depending on the parameter $\vartheta$, we shall have three different regimes of interactions, which lead to three different asymptotic limits of $Z_t^\eps$ as $\eps\to0$, i.e.,
\begin{equation}\label{regime}
\left\{\begin{aligned}
&\text{Regime 1: } \,\, \vartheta<1/2;\\
&\text{Regime 2: } \,\, \vartheta>1/2;\\
&\text{Regime 3: } \,\, \vartheta=1/2.
\end{aligned}\right.
\end{equation}
In Regime 1, the time average of the drift will be the dominant term, and  no homogenization but an
explicit deterministic contribution appears in the limit dynamics;
whereas in Regime 2, the space average of the drift will be the leading term and
homogenization   will take place, the limit of $Z_t^\eps$ will be an Ornstein-Uhlenbeck type process with homogenized diffusion coefficient given in terms of the solution of an auxiliary infinite-horizon Kolmogorov equation;
finally, in Regime 3 both the deterministic contribution and
  homogenization will occur together, see {\bf Theorem \ref{main2}} below.

\vspace{1mm}
Compared with previous works relying on time discretization, our method for proving the above convergence is quite unified, and the assumptions on the coefficients are very weak. In particular, periodic, quasi-periodic and almost periodic coefficients are allowed, see Remark \ref{th1} and Remark \ref{r2} below. Furthermore, explicit rates of the above convergence   are also obtained as straightforward byproducts of our arguments, which do not depend on the   regularities of the coefficients with respect to the time component and the fast variable. The main tool we shall use is the following
infinite-horizon time-inhomogeneous Kolmogorov equation in $[0,\infty)\times\mR^{d_1}$ and with a parameter:
\begin{align}\label{pde0}
\partial_t u(t,x,y)+\sL_0u(t,x,y)=-f(t,x,y),
\end{align}
where $y\in\mR^{d_2}$ is a parameter, and the operator $\sL_0$ is given by
\begin{align}\label{L0}
\sL_0:=\sL_0(t,x,y)
=b(t,x,y)\cdot\nabla_x+\frac{1}{2}\tr\left(\sigma\sigma^*(t,x,y)\cdot\nabla_x^2\right),
\end{align}
which is the  generator of  $X_{s,t}^y(x)$ satisfying the frozen equation (\ref{froz}). We point out that such kind of parabolic equation has been used by Bensoussan, Lions and Papanicolaou  to study
the asymptotic limit for the diffusion
\begin{align}\label{L00}
\dif X_t^\eps=\frac{1}{\eps}b(t/\eps^2,X_t^\eps/\eps)\dif t+\frac{1}{\sqrt\eps}\sigma(t/\eps^2,X_t^\eps/\eps)\dif W_t
\end{align}
when the coefficients $\sigma$ and $b$ are time-space periodic, see \cite[Chapter 3]{AJG}. In the homogenization literature, it was also used   as a corrector when the coefficients are space-time stationary
(e.g., almost periodic), see the excellent paper \cite{PAB}, see also \cite{ABM,BS,FB} and the references therein. Note that we are studying the equation (\ref{pde0}) in the whole space but not on compact sets, and there are no initial or terminal conditions. To be well-defined, we assume  a  confining drift (see  $(\mathbf{ H_1})$ below) to guarantee the existence of an evolution
system of invariant measures $\{\mu_t^y\}_{t\geq 0}$ for the corresponding process $X_{s,t}^y(x)$, and a necessary condition  for the existence of a solution $u$ is centering for the right hand side $f$:
\begin{align}\label{cen}
\int_{\mR^{d_1}}f(t,x,y)\mu_t^y(\dif x)=0,\quad\forall (t,y)\in\mR_+\times\mR^{d_2}.
\end{align}
 We shall establish the well-posedness and the regularities of the solution $u$ with respect to the $x$ variable as well as the parameter $y$, see  Theorem \ref{cauchy1}.
Then we use $u(t,x,y)$ to serve as a corrector to  establish certain weak and strong fluctuation estimates of double averaging type and functional central
limit theorem type for the non-autonomous system, see Lemmas \ref{xx}, \ref{xxx} and  \ref{weaf}, respectively.
 We believe that our result for the inhomogeneous equation (\ref{pde0}) can also be used to treat similar topics pertaining to
(\ref{poi}) in the non-autonomous context (e.g., numerical approximation  of invariant measures of inhomogeneous SDEs as in \cite{NX,MST,SX}), and   to study the asymptotic limit of the system (\ref{L00}) in the whole space without the space-periodic assumption, thus the result of Theorem \ref{cauchy1} is of independent interest.

\vspace{1mm}
The rest of this paper is organized as follows. Section 2 states the main assumptions and results.
Section 3 is devoted to studying the time-inhomogeneous Kolmogorov equation. In Section 4, we establish some fluctuation
lemmas  and prove the
weak and strong convergence in the averaging principle. Finally, we  give the proof of the functional central limit type theorem in Section 5.

\vspace{2mm}
\noindent{\bf Notations:} Throughout this paper, we use the following convention:   $C$ with or without subscripts will denote positive constants, whose values may change in different places, and whose dependence on parameters can be traced from the calculations.
Given two functions $f(t)$ and $g(t)$, we say that $f(t) \propto g(t)$ as $t\to\infty$ if there exists a non-zero constant $C$ such that  $\lim_{t\to\infty}f(t)/g(t)=C$. Finally, we introduce the following spaces of functions used in this paper.

\begin{itemize}
\item For a function $f(t,s,x ,y,z)$ defined on $\mR_+^2\times \mR^{d_1}\times \mR^{2d_2},$ by $f\in L_p^\infty:=L_p^\infty(\mR_+^2\times \mR^{d_1}\times \mR^{2d_2})$ we mean that there exist constants $C,p>0$   such that
$$|f(t,s,x,y,z)|\leq C(1+|x|^p),\;\forall t,s\geq 0,x\in\mR^{d_1},y,z\in\mR^{d_2}.$$
\item For $0<\alpha,\beta\leq 1$, the space $C_p^{\alpha,\beta}(\mR^{d_1}\times \mR^{d_2})$ consists of all functions $f\in L_p^\infty(\mR^{d_1}\times \mR^{d_2})$ that are $\alpha$-local H\"older continuous with polynomial growth in $x$ and  $\beta$-H\"older continuous   in $y$, i.e.,  there exist constants $C,p>0$ such that for any $x_1,x_2\in\mR^{d_1}$ and $y_1,y_2\in\mR^{d_2}$,
	\begin{align*}
		|f(x_1,y_1)-f(x_2,y_2)|&\leq C\left[\big(|x_1-x_2|^\alpha\wedge1\big)+\big(|y_1-y_2|^\beta\wedge1\big)\right]\\
		&\qquad\qquad\times\big(1+|x_1|^p+|x_2|^p\big).
	\end{align*}
When  $\alpha,\beta\geq1$, we denote
$$
C_p^{\alpha,\beta}(\mR^{d_1}\times \mR^{d_2}):=\Big\{f: \partial_x ^{[\alpha]}\partial_y^{[\beta]}f\in C_p^{\alpha-[\alpha],\beta-[\beta]}(\mR^{d_1}\times \mR^{d_2}) \Big\}.
$$

\item For $0<\delta,\vartheta\leq 1$, the space $C_p^{\delta,\vartheta,\alpha,\beta}(\mR_+^2\times\mR^{d_1}\times \mR^{d_2})$ contains all functions $f$ such that for every fixed $t,s$, $f(t,s,\cdot,\cdot)\in C_p^{\alpha,\beta}(\mR^{d_1}\times \mR^{d_2})$, and for every fixed $x,y$,
$f(\cdot,\cdot,x,y)\in C_b^{\delta,\vartheta}(\mR_+^2)$ (the subscript $b$ stands for boundedness).

 Similarly for $\gamma>0,$ the space
$C_p^{\delta,\vartheta,\alpha,\beta,\gamma}(\mR_+^2\times\mR^{d_1}\times \mR^{2d_2})$ contains all functions $f$ such that for every fixed $z\in\mR^{d_2}$, $f(\cdot,\cdot,\cdot,\cdot,z)\in C_p^{\delta,\vartheta,\alpha,\beta}(\mR_+^2\times\mR^{d_1}\times \mR^{d_2})$, and for every fixed $(t,s,x,y)$, $f(t,s,x,y,\cdot)\in C_b^{\gamma}(\mR^{d_2})$.

\item
	For $0<\alpha,\beta\leq 1$, the space $\mathbf{C}_p^{\alpha,\beta}(\mR^{d_1}\times \mR^{d_2})$ consists of all functions $f$ such that for some
	$C, p>0$,
	\begin{align*}
	|f(x_1,y_1)-f(x_2,y_2)|&\leq C\left[\big(|x_1-x_2|^\alpha\wedge1\big)+\big(|y_1-y_2|^\beta\wedge1\big)\right]\\
	&\,\, \,\,\times\big(1+|x_1|^p+|x_2|^p+|y_1|^p+|y_2|^p\big), \forall x_1,x_2\in\mR^{d_1}, y_1,y_2\in\mR^{d_2},
	\end{align*}
For $\vartheta>0$, the space $\mathbf{C}_{p}^{\vartheta,\alpha,\beta}(\mR_+\times\mR^{d_1}\times \mR^{d_2})$ contains all functions $f$ such that for every fixed $t$, $f(t,\cdot,\cdot)\in \mathbf{C}_p^{\alpha,\beta}(\mR^{d_1}\times \mR^{d_2})$, for every fixed $x,y$,
		$f(\cdot,x,y)\in C_b^{\vartheta}(\mR_+)$.
	
	\end{itemize}

\section{Assumptions and main results}

We study the asymptotic behavior of the multi-scale system (\ref{sde0}) in the whole space $\mR^{d_1+d_2}$.
While we
expect that the results of this paper would still hold if $b$ and $\sigma$ were assumed to be periodic in the $x$-variable, we instead consider the case when the drift has a
confining effect, leading to  ergodicity of the frozen process.
To be precise, we shall make the following assumptions:

\begin{enumerate}[$(\mathbf{ H_1})$]
\item The coefficient $\sigma(t,x,y)$ is  locally elliptic in $x$ uniformly with respect to $(t,y)$, i.e., there exist $\lambda, p>1$ such that for any $(t,y)\in\mR_+\times\mR^{d_2},$
$$
\lambda^{-1}(1+|x|)^{-p}|\xi|^{2} \leq |\sigma(t,x,y) \xi|^2 \leq \lambda(1+|x|)^{p} |\xi|^{2}, \quad \forall \,x,\xi \in \mathbb{R}^{d_1},
$$
and there exist $c_0,c_1> 0$ such that
\begin{align}\label{ex6}
|\sigma(t,x,y)|^2_{HS}+\<x,b(t,x,y)\>\leq -c_0|x|^2+c_1.
\end{align}
\end{enumerate}

Under $(\mathbf{ H_1})$,  it was shown that the time-inhomogeneous SDE (\ref{froz}) admits an evolution system of invariant measures $\{\mu_t^y\}_{t\geq 0}$, see e.g. \cite{DR1} and \cite[Theorem 5.4]{KLL}. Concerning the coefficients of the slow equation, we assume:

\begin{enumerate}[$(\mathbf{ H_2})$]
\item The coefficient $G(t,x,y)$ is bounded and elliptic in $y$ uniformly with respect to $(t,x)$, i.e., there exists $\lambda>1$ such that for any $(t,x)\in\mR_+\times\mR^{d_1}$,
$$
\lambda^{-1}|\xi|^2\leq |G(t,x,y)\xi|^2\leq\lambda|\xi|^2,\ \ \forall\,\xi,y\in\mR^{d_2}.
$$
Moreover, there exist  functions $\kappa_i:\mR_+\to\mR_+$, $i=1,2$, satisfying $\kappa_i(T)\to 0$ as $T\to \infty$, $\bar {\bar F}:\mR^{d_2}\to\mR^{d_2}$ and $\bar {\bar G}:\mR^{d_2}\to\mR^{d_2\times d_2}$ such that for all $T>0$ and $y\in\mR^{d_2}$,
\begin{align}\label{bF}
\left|\frac{1}{T}\int_{0}^{T}\!\!\int_{\mR^{d_1}}F(t,x,y)\mu_t^y(\dif x)\dif t-\bar{\bar F}(y)\right|\leq \kappa_1(T),
\end{align}
and
\begin{align}\label{bGG}
\left|\frac{1}{T}\int_{0}^{T}\!\!\int_{\mR^{d_1}}GG^*(t,x,y)\mu_t^y(\dif x)\dif t- \bar{\bar G}\bar{\bar G}^*(y)\right|_{HS}\leq \kappa_2(T),
\end{align}
where $|\cdot|_{HS}$ is the Hilbert-Schmidt norm of the matrix.
\end{enumerate}

We shall show that as $\eps\to0$, the slow process $Y_t^\eps$ of the system (\ref{sde0}) will converge weakly  to $\bar Y_t$ which satisfies the following double averaged equation:
\begin{align}\label{wbsde}
\dif \bar Y_t=\bar{\bar F}(\bar Y_t)\dif t + \bar {\bar G}(\bar Y_t)\dif W_t^2,\quad\bar Y_0=y\in\mR^{d_2},
\end{align}
where $\bar{\bar F}$ and  $ \bar {\bar G}$ are given in (\ref{bF}) and (\ref{bGG}), respectively, i.e.,
\begin{align*}
\bar{\bar F}(y)=\lim\limits_{T\to\infty}\frac{1}{T}\int_0^{T}\!\!\int_{\mR^{d_1}}F(t,x,y)\mu_t^y(\dif x)\dif t,
\end{align*}
and
$$
\bar{\bar G}\bar{\bar G}^*(y)=\lim\limits_{T\to \infty}\frac{1}{T}\int_0^{T}\int_{\mR^{d_1}}GG^*(t,x,y)\mu_t^y(\dif x)\dif t.$$
The following is the first main result of this paper.

\bt\label{main1}
Let {\bf $(\mathbf{ H_1})$} and {\bf $(\mathbf{ H_2})$} hold. Assume that $\sigma, b, F, G\in C_p^{\alpha/2,\alpha,\beta}$ with $0<\alpha,\beta$.
Then for any $T>0$ and $\varphi\in C_b^{2+\b}$, we have
\begin{align}\label{weak}
\sup_{t\in[0,T]}|\mE\varphi(Y_t^\eps)-\mE\varphi(\bar Y_t)|\leq C_{T,\varphi}\bigg(\eps^{(\beta /2)\wedge1}+\sup_{\substack{1\leq i\leq 2\\t\in [0,T]}}\Big[t\cdot\kappa_i(t/\eps)\Big]\bigg),
\end{align}
where $\bar Y_t$ satisfies SDE (\ref{wbsde}), the functions $\kappa_1$ and $\kappa_2$ are given in (\ref{bF}) and (\ref{bGG}), respectively, and $C_{T,\varphi}>0$ is a constant independent of $\a$ and $\eps$.
\et
We give the following remark for the above result.

\br\label{th1}
(i)
The assumptions   (\ref{bF}) and (\ref{bGG})    hold when the coefficients are periodic, quasi-periodic or
almost periodic in time $t$ uniformly with respect to $x,y$. Let us explain this for (\ref{bF}) in the almost periodic case, since the periodic and quasi-periodic functions are also almost periodic.  To be precise, since $b,\sigma$ are uniformly almost periodic, similarly  as in \cite[Theorem 6.3]{Ce} we have that $\mu_t^y$ is almost periodic for any fixed $y\in\mR^{d_2}$. Let
\begin{align}\label{btF}
\bar F(t,y):=\int_{\mR^{d_1}}F(t,x,y)\mu_t^y(\dif x).
\end{align}
By  similar arguments as  in \cite[Lemma 7.1]{Ce} and the fact that $F(t,x,y)$ is uniformly almost periodic, we know that $\bar F(t,y)$ is almost periodic in $t$ uniformly with respect to $y$. Thus, (\ref{bF}) holds with
\begin{align*}
\bar{\bar F}(y):=\lim\limits_{T\to\infty}\frac{1}{T}\int_0^{T}\bar F(t,y)\dif t,
\end{align*}
see e.g. \cite[Theorem 6.11]{Co}.

(ii) If the coefficients of system (\ref{sde0})  are periodic functions, then the conditions (\ref{bF}) and  (\ref{bGG})   hold with
\begin{align}\label{per}
\kappa_1(T)\propto 1/T\quad\text{and}\quad\kappa_2(T)\propto 1/T
\end{align}
as $T\to\infty$, see e.g. \cite[(1.9')]{Il}.
Thus, the estimate (\ref{weak}) becomes
\begin{align*}
\sup_{t\in[0,T]}|\mE\varphi(Y_t^\eps)-\mE\varphi(\bar Y_t)|\leq C_{T,\varphi}\,\eps^{(\beta /2)\wedge1},
\end{align*}
which coincides with the result for classical multi-scale autonomous stochastic systems.
\er

To establish the strong convergence in the averaging principle, we need to further assume that the diffusion coefficient in the slow equation does not depend on the fast motion. We make the following assumption:
\begin{enumerate}[$(\mathbf{\hat H_2})$]
\item The coefficient
$G(t,x,y)\equiv G(t,y)$ (does not depend on $x$) is bounded and elliptic in $y$ uniformly with respect to $t$, the drift $F$ satisfies (\ref{bF})
and there exist   functions $\hat\kappa_2:\mR_+\to\mR_+$ satisfying $\hat\kappa_2(T)\to 0$ as $T\to \infty$ and   $\bar G:\mR^{d_2}\to\mR^{d_2\times d_2}$ such that for all $T>0$ and $y\in\mR^{d_2}$,
\begin{align}\label{bG}\frac{1}{T}\int_{0}^{T}\left|G(t,y)-{\bar G}(y)\right|_{HS}^2\dif t\leq \hat\kappa_2^2(T),\end{align}
where $|\cdot|_{HS}$ is the Hilbert-Schmidt norm.
\end{enumerate}

We prove that $Y_t^\eps$ will converge strongly (in $L^q(\Omega)$-sense) as $\eps\to0$ to $\bar Y_t$ which satisfies
\begin{align}\label{bsde}
\dif \bar Y_t=\bar{\bar F}(\bar Y_t)\dif t+ \bar G(\bar Y_t)\dif W_t^2,\quad\bar Y_0=y\in\mR^{d_2},
\end{align}
where   $\bar{\bar F}$
and $\bar G$ are given in (\ref{bF}) and (\ref{bG}), respectively. We have the following result.

\bt\label{main11}
Let {\bf $(\mathbf{ H_1})$} and {\bf $(\mathbf{\hat H_2})$} hold. Assume that $\sigma, b, F\in C_p^{\alpha/2,\alpha,\beta}$ and $G\in C_b^{\alpha/2,1}$ with $0<\a,\b$.
Then for any $T>0$ and $q\geq 1$,  we have
\begin{align}\label{strong}
\sup_{t\in[0,T]}\mE|Y_t^\eps-\bar Y_t|^q\leq C_{T}\Big(\eps^{(\beta\wedge1) q/2}+\sup_{t\in[0,T]}\big[t\cdot\kappa_1(t/\eps)\big]^q+\sup_{t\in[0,T]}
\big[t\cdot\hat\kappa_2(t/\eps)\big]^q\Big),\end{align}
where $\bar Y_t$ satisfies SDE (\ref{bsde}), the functions $\kappa_1$ and $\hat\kappa_2$ are given in (\ref{bF}) and (\ref{bG}), respectively, and $C_{T}>0$ is a constant independent of $\a$ and $\eps$.
\et
\br
(i) As in (\ref{gg}), in general the independence of $G$ with respect to the $x$-variable is necessary for the strong convergence to hold. In particular, if $G$ is  asymptotically constant in time, then (\ref{bG}) holds. The main reason for assuming (\ref{bG}), compared with the assumption (\ref{bGG}) for the weak convergence,  is that we need to further control  the convergence of the martingale part (see the third term in (\ref{mar})), while the expectation of the martingale part equals  zero in proving the weak convergence (see (\ref{cI}) below).

(ii)  Let us give some comments for the rate of convergence in (\ref{strong}). For simplicity, we consider $G(t,x,y)\equiv G(y)$ and thus $\hat\kappa_2=0$ in (\ref{bG}).
The strong convergence rate of $Y_t^\eps$ to $\bar Y_t$ contains   two parts: $\eps^{\beta/2}$ controls the difference between $F(t/\eps, X_t^\eps, Y_t^\eps)$ and its space average $\bar F(t/\eps, Y_t^\eps)$ defined by (\ref{btF}), while $\kappa_1(t/\eps)$ controls the difference between $\bar F(t/\eps, Y_t^\eps)$ and its time average $\bar{\bar F}(Y_t^\eps)$, see Lemmas \ref{xx} and \ref{xxx}. This will be important for determining the normalization constant below.
\er

By a localization technique, we can also drop the uniform ellipticity and the boundedness assumptions on the coefficients with respect to the slow variable.  We assume:	
\begin{enumerate}[$(\mathbf{\hat H_{2;loc}})$]
	\item The coefficient
	$G(t,x,y)\equiv G(t,y)$ (does not depend on $x$) is locally bounded and elliptic in $y$ uniformly with respect to $t$,
	 i.e., there exist constants $\lambda, p>1$ such that for any $(t,y)\in\mR_+\times\mR^{d_2},$
	$$
	\lambda^{-1}(1+|y|)^{-p} |\xi|^{2} \leq |G(t,y) \xi|^2 \leq 	\lambda (1+|y|)^{p}|\xi|^{2}, \quad \forall \xi \in \mathbb{R}^{d_2}.
	$$
Moreover, there exist $p>1$, functions $\hat\kappa_i:\mR_+\to\mR_+$, $i=1,2$ satisfying $\hat\kappa_i(T)\to 0$ as $T\to \infty$, $\bar {\bar F}:\mR^{d_2}\to\mR^{d_2}$ and $ \bar G:\mR^{d_2}\to\mR^{d_2\times d_2}$ such that for all $T>0$ and $y\in\mR^{d_2}$,
\begin{align}\label{bF0}
\left|\frac{1}{T}\int_{0}^{T}\!\!\int_{\mR^{d_1}}F(t,x,y)\mu_t^y(\dif x)\dif t-\bar{\bar F}(y)\right|\leq \hat\kappa_1(T)(1+|y|^p),
\end{align}
and
\begin{align}\label{bG0}
\frac{1}{T}\int_{0}^{T}\left|G(t,y)-{\bar G}(y)\right|_{HS}^2\dif t\leq \hat\kappa_2^2(T)(1+|y|^p),
\end{align}
where $|\cdot|_{HS}$ is the Hilbert-Schmidt norm.
\end{enumerate}	

We have the following result.

\bc\label{coro}
Let {\bf $(\mathbf{ H_{1}})$} and {\bf $(\mathbf{\hat H_{2;loc}})$} hold. Assume that $\sigma, b\in C_{p}^{\alpha/2,\alpha,\beta}$, $F\in \mathbf{C}_{p}^{\alpha/2,\alpha,\beta}$ and $G\in \mathbf{C}_{p}^{\alpha/2,1}$ with $0<\a,\b$, and that the following moment estimate holds:

{\bf($\Theta$)} For any $T>0$, there exist constants $\theta \geqslant 1$ and $C_{\theta,T}>0$ such that
$$
\sup _{\varepsilon \in(0,1)} \mathbb{E}\left(\sup _{t \in[0, T]}\left(\left|Y_t^{\varepsilon}\right|^\theta+\left|\bar{Y}_t\right|^\theta\right)\right) \leqslant C_{\theta,T}<\infty .
$$
Then for every $q<\theta$, we have
$$
\lim _{\varepsilon \rightarrow 0} \sup _{t \in[0, T]} \mathbb{E}\left|Y_t^{\varepsilon}-\bar{Y}_t\right|^q=0,
$$
where $\bar{Y}_t$ satisfies SDE \eqref{bsde}.
\ec

\br
The advantage of Corollary \ref{coro} is that we only need to show the a priori moment estimate {\bf($\Theta$)} in order to guarantee  the strong convergence in the averaging principle for SDE (\ref{sde0}) with locally H\"older continuous drifts. For example, if $G$ and $F$ satisfy
$$
|G(t,y)|_{HS}\leq C_0(1+|y|)\quad {\text{and}}\quad \<F(t,x,y),y\>\leq C_0(1+|x|^p+|y|^2),
$$
where $p>1$ and $C_0>0$ is independent of $t$, then {\bf($\Theta$)} follows directly by It\^o's formula and the moments estimate of the fast process $X_t^\eps$ (which is guaranteed by (\ref{ex6}) since the constant is independent of $y$).
\er

Next, we proceed to study the small
fluctuations of $Y_t^\eps$ from its average $\bar Y_t$. For simplicity, we assume that:

\begin{enumerate}[${\bf ({\tilde H_2})}$]
\item The coefficient $G(t,x,y)\equiv G(y)$ (does not depend on $(t,x)$) is bounded and elliptic,   and there exists a $0<\vartheta\leq 1$ such that for every $t>0$,
$$t\cdot\kappa_1(t/\eps)
\propto\eps^\vartheta$$
as $\eps\to0$,
 where $\kappa_1$ is given in (\ref{bF}).
\end{enumerate}

Under assumption {$\bf (\tilde{H}_2$)} and when $\beta=1$ in Theorem \ref{main1}, we know that $Y_t^\eps$ converges strongly to $\bar Y_t$ with a  rate $\sqrt{\eps}+\eps^{\vartheta}$. Recall that the normalized deviation process $Z_t^\eps$ is defined by (\ref{zep}) and we consider Regime 1-3 in (\ref{regime}). To characterize  the limit of the difference between $\bar F$ and its time average $\bar{\bar F}$, we make the following assumption:

\begin{enumerate}[${\bf ({\tilde H_3})}$]
\item There exist functions
$\tilde\kappa_1:\mR_+\to\mR_+$ satisfying
$\tilde\kappa_1(T)\to0$ as $T\to\infty$ and $\Upsilon_\vartheta\in C_b^{1+\beta}(\mR^{d_2})$  such that for all $T>0$ and $y\in\mR^{d_2}$,
\begin{align}\label{U1}
&\left|
\frac{1}{T}\int_0^T
\left[
\frac{\bar F(t/\eps,y)-\bar{\bar F}(y)}
{\eps^\vartheta}
-\Upsilon_\vartheta(y)
\right]\dif t
\right|\leq
\tilde\kappa_1(T/\eps),
\end{align}
where $\bar F$ is defined by (\ref{btF}).
\end{enumerate}

In order to characterize the homogenized limit of the space average of the drift, we need to consider the following time-inhomogeneous Kolmogorov equation in $\mR_+\times\mR^{d_1}$:
\begin{align}\label{PPhi}
\partial_t\Phi(t,x,y)+\sL_0(t,x,y)\Phi(t,x,y)=-[F(t,x,y)-\bar F(t,y)],
\end{align}
where $y\in\mR^{d_2}$ is regarded as  a parameter, and the operator $\sL_0(t,x,y)$ is given by (\ref{L0}). Note that the right hand side satisfies the centering condition (\ref{cen}).
Under our assumptions below, there exists a unique solution $\Phi$ to equation (\ref{PPhi}). Define
\begin{align}\label{Sigma}
\frac{1}{2}\bar\Sigma\bar\Sigma^*(t,y):=\int_{\mR^{d_1}}[F(t,x,y)-\bar F(t,y)]\Phi^*(t,x,y)\mu_t^y(\dif x).
\end{align}
We further make the following  assumption.

\vspace{2mm}
\begin{enumerate}[$ {\bf ({\tilde H_4})} $]
\item There exist functions $\tilde\kappa_2:\mR_+\to\mR_+$ satisfying $\tilde\kappa_2(T)\to 0$ as $T\to \infty$ and  $\bar{\bar\Sigma}:\mR^{d_2}\to\mR^{d_2\times d_2}$
such that for all $T>0$ and $y\in\mR^{d_2}$,
\begin{align}\label{bsigma}\left|\frac{1}{T}\int_{0}^{T}
{\bar\Sigma}{\bar\Sigma}^*(t,y)
\dif t-\bar{\bar\Sigma}\bar{\bar\Sigma}^*(y)\right|_{HS}\leq \tilde\kappa_2(T).\end{align}
\end{enumerate}
In particular,
\begin{align} \label{bsigma'}
\bar{\bar\Sigma}\bar{\bar\Sigma}^*(y)=\lim\limits_{T\to \infty}\frac{1}{T}\int_0^{T}\bar\Sigma\bar\Sigma^*(t,y)\dif t.
\end{align}
Corresponding to Regime $k$ ($k=1,2,3$) in (\ref{regime}), we shall show that
 $Z_t^\eps$ converges weakly to $\bar Z_t^k$ ($k=1,2,3$), where
\begin{align}
&\dif \bar Z^1_{t}=\nabla_y\bar{\bar F}(\bar Y_t)\cdot\bar Z^1_{t}\dif t+\Upsilon_\vartheta(\bar Y_t)\dif t+\nabla_y G(\bar Y_t)\cdot\bar Z^1_{t}\dif W_t^2,\label{bze1}\\
&\dif \bar Z^2_{t}=\nabla_y\bar{\bar F}(\bar Y_t)\cdot\bar Z^2_{t}\dif t+\nabla_y G(\bar Y_t)\cdot\bar Z^2_{t}\dif W_t^2+\bar{\bar\Sigma}(\bar Y_t)\dif \tilde W_t,\label{bze2}\\
&\dif \bar Z^3_{t}=\nabla_y\bar{\bar F}(\bar Y_t)\cdot\bar Z^3_{t}\dif t+\Upsilon_{1/2}(\bar Y_t)\dif t+\nabla_y G(\bar Y_t)\cdot\bar Z^3_{t}\dif W_t^2+\bar{\bar\Sigma}(\bar Y_t)\dif \tilde W_t,\label{bze3}
\end{align}
respectively, where $\tilde W_t$ is another Brownian motion independent of $W_t^2$, $\bar{\bar F}(y)$ is defined by (\ref{bF}),   and $\Upsilon_\vartheta(y)$ and $\bar{\bar \Sigma}(y)$ are  given in  (\ref{U1}) and (\ref{bsigma'}).
We have the following result. In particular, the periodic setting  belongs to Regime 2, see   Remark \ref{r2} below.

\bt\label{main2}
Let {\bf $(\mathbf{ H_1})$} and {\bf $(\mathbf{ \tilde H_2})$}   hold,  $T>0$,  $\sigma, b, F\in C_p^{\alpha/2,\alpha,1+\beta}$  and $G\in C_b^{1+\beta}$ with  $\alpha,\b\in (0,1]$.

\vspace{1mm}
\noindent
(i) (Regime 1: $\vartheta\in(0,1/2)$) Assume that $(\mathbf{ \tilde H_3})$ holds.
Then for every $\varphi\in C_b^4(\mR^{d_2})$, we have
$$
\sup_{t\in[0,T]}|\mE\varphi(Z_t^\eps)-\mE\varphi(\bar Z^1_{t})|\leq C_{T,\varphi}
\Big(\eps^{\vartheta\wedge(1-2\vartheta)\wedge\beta/2}+\sup\limits_{t\in[0,T]}t\cdot\tilde\kappa_1(t/\eps)\Big);
$$
(ii) (Regime 2: $\vartheta\in(1/2,1]$) Assume that $(\mathbf{ \tilde H_4})$ holds.
Then for every $\varphi\in C_b^4(\mR^{d_2})$, we have
$$
\sup_{t\in[0,T]}|\mE\varphi(Z_t^\eps)-\mE\varphi(\bar Z^2_{t})|\leq C_{T,\varphi}\Big(\eps^{(\vartheta-1/2)\wedge\beta /2}+\sup\limits_{t\in[0,T]}t\cdot\tilde\kappa_2(t/\eps)\Big);
$$
(iii) (Regime 3: $\vartheta=1/2$) Assume that $(\mathbf{ \tilde H_3})$ with $\vartheta=1/2$ and $(\mathbf{ \tilde H_4})$ hold.
Then for every $\varphi\in C_b^4(\mR^{d_2})$, we have
$$
\sup_{t\in[0,T]}|\mE\varphi(Z_t^\eps)-\mE\varphi(\bar Z^3_{t})|\leq C_{T,\varphi}\bigg(\eps^{\beta/2}
+\sup\limits_{\substack{1\leq i\leq 2\\t\in[0,T]}}t\cdot\tilde\kappa_i(t/\eps)\bigg);
$$
where $C_{T,\varphi}>0$ is  a constant independent of $\a$ and $\eps$, 
and $\tilde\kappa_1,\tilde\kappa_2$ are given in (\ref{U1})  and (\ref{bsigma}), respectively.
\et

\br\label{r2}
(i) Notably, when $G(y)\equiv\mI_{d_2}$ (the identity matrix), the limit process switches from a deterministic correction (Regime 1) to an Ornstein-Uhlenbeck type process with homogenized noise (Regime 2) when $\vartheta$ crosses $1/2$.

(ii)
In view of (\ref{per}), the case when all the coefficients are periodic in time belongs to Regime 2. In this situation, $(\mathbf{ \tilde H_4})$ holds with
\begin{align}\label{88}
t\cdot\tilde\kappa_2(t/\eps)
\propto t\cdot\eps/t=\eps
\end{align}
as $\eps\to0$, and  thus we obtain
$$
\sup_{t\in[0,T]}|\mE\varphi(Z_t^\eps)-\mE\varphi(\bar Z^2_{t})|\leq C_{T,\varphi}\,\eps^{\beta/2}.
$$
To verify (\ref{88}) and according to  \cite[(1.9')]{Il}, we only need to show that the function ${\bar\Sigma}{\bar\Sigma}^*(t,y)$ defined by (\ref{Sigma}) is periodic.
 Assume that $b,\sigma$ are $\tau$-periodic functions, by \cite[Proposition 2.10]{LLZ}, we have that $\mu_t^y$ is also $\tau$-periodic. Thus, by the definition (\ref{btF}) and the periodicity of $F$, we have that $\bar F$ is $\tau$-periodic.
It remains to prove the periodicity of  $\Phi(t,x,y)$. Recall that for all $t\in\mR_+,$
\begin{align*}
\partial_t\Phi(t,x,y)+\sL_0(t,x,y)\Phi(t,x,y)=-[F(t,x,y)-\bar F(t,y)].
\end{align*}
Thus,
\begin{align*}
\partial_t\Phi(t+\tau,x,y)+\sL_0(t+\tau,x,y)\Phi(t+\tau,x,y)=-[F(t+\tau,x,y)-\bar F(t+\tau,y)].
\end{align*}
Since $F,b,\sigma$ are $\tau$-periodic, we have $\sL_0(t+\tau,x,y)=\sL_0(t,x,y), \bar F(t+\tau,y)=\bar F(t,y).$ Hence we have
\begin{align*}
\partial_t\Phi(t+\tau,x,y)+\sL_0(t,x,y)\Phi(t+\tau,x,y)=-[F(t,x,y)-\bar F(t,y)].
\end{align*}
By the uniqueness of the solution of the Poisson equation in Theorem \ref{cauchy1}, we obtain $\Phi(t+\tau,x,y)=\Phi(t,x,y).$

\er

\section{Infinite-horizon  time-inhomogeneous Kolmogorov equation}

The aim of this section is to study the well-posedness and the regularities of the solution of the  time-inhomogeneous PDE (\ref{pde0}).
Recall that $\sL_0$ is the infinitesimal generator of the frozen  process $X_{s,t}^y(x)$ which satisfies the SDE (\ref{froz}). Under the assumption {\bf $(\mathbf{ H_1})$}, there exists an evolution system of invariant measures $\mu^y_t$ such that  for every $y\in\mR^{d_2}$, $t\geq s>0$ and $\varphi\in C_p(\mR_+\times\mR^{d_1})$,
\begin{align}\label{ergo2}
\int_{\mR^{d_1}}P^y_{s,t}\varphi(t,x)\mu_s^y(\dif x)=\int_{\mR^{d_1}}\varphi(t,x)\mu_t^y(\dif x),
\end{align}
and there exist constants $C,p,\delta>0$ independent of $(t,y)$ such that
\begin{align}\label{decay2}
\left|P^y_{s,t}\varphi(t,x)-\int_{\mR^{d_1}}\varphi(t,z)\mu_t^y(\dif z)\right|\leq C(1+|x|^p)e^{-\delta(t-s)},
\end{align}
where
\begin{align}\label{semi}
P^y_{s,t}\varphi(t,x):=\mE\varphi\big(t,X_{s,t}^y(x)\big).
\end{align}
see e.g. \cite{DR1} and \cite[Theorem 5.4]{KLL}. In addition, it can be shown that $\{\mu_t^y\}_{t\geq 0}$ satisfy  the Fokker-Planck equation that
\begin{align}\label{FP}
(\p_t+\sL_0)^*\mu_t^y=0,
\end{align}
where $(\p_t+\sL_0)^*$ is the formal adjoint operator.
In fact, for every $\varphi\in C^{1,2}_p(\mR_+\times\mR^{d_1})$, using It\^o's formula and taking expectation, we have
\begin{align*}
\mE\varphi\big(t,X_{s,t}^y(x)\big)=\varphi(s,x)
+\mE\left(\int_{s}^t\p_r\varphi\big(r,X_{s,r}^y(x)\big)+\sL_0\varphi\big(r,X_{s,r}^y(x)\big)\dif r\right).
\end{align*}
Taking integral with respect to $\mu_s^y(\dif x)$ from both sides of the above equality and by (\ref{ergo2}), we deduce that
\begin{align*}
&\int_{\mR^{d_1}}\varphi(t,x)\mu_t^y(\dif x)=\int_{\mR^{d_1}}\!P^y_{s,t}\varphi(t,x)\mu_s^y(\dif x)\\
&=\int_{\mR^{d_1}}\!\varphi(s,x)\mu_s^y(\dif x)+\int_{s}^t\!\!\int_{\mR^{d_1}}P_{s,r}^y\big(\p_r\varphi+\sL_0\varphi\big)(r,x)\mu_s^y(\dif x)\dif r\\
&=\int_{\mR^{d_1}}\!\varphi(s,x)\mu_s^y(\dif x)+\int_{s}^t\!\!\int_{\mR^{d_1}}\big(\p_r\varphi+\sL_0\varphi\big)(r,x)\mu_r^y(\dif x)\dif r,
\end{align*}
which in turn implies (\ref{FP}). As a result, taking integral with respect to $\mu^y_t$ for (\ref{pde0}), we see that
the function $f(t,x,y)$ needs to
satisfy the centering condition (\ref{cen}).

We have the following result.

\bt\label{cauchy1}
Assume that $\sigma, b\in C_p^{\alpha/2,\alpha,\beta}$ with $0<\alpha\leq 1$ and $\beta\geq 0$. Then, for every $f\in C_p^{\alpha/2,\alpha,\beta}$ satisfying (\ref{cen}),
there exists a unique solution $u\in C_p^{1+\alpha/2,2+\alpha,\beta}$ to the equation  (\ref{pde0}) which also satisfies the centering condition  (\ref{cen}). Moreover, $u$ admits the probability representation that
\begin{align}\label{solv}
u(t,x,y)=\int_t^\infty P^y_{t,r}f(r,x,y)\dif r,\quad \forall(t,x,y)\in\mR_+\times\mR^{d_1+d_2},
\end{align}
where $P^y_{t,r}$ is defined by (\ref{semi}).
\et

\begin{proof}
Let us first show that $u$ in (\ref{solv}) is well-defined and satisfies (\ref{cen}).
Since $f$ satisfies the centering condition (\ref{cen}), we have by (\ref{decay2}) that for every $r>t$,
\begin{align*}
|P^y_{t,r}f(r,x,y)|&=\bigg|P^y_{t,r}f(r,x,y)-\int_{\mR^{d_1}}f(r,z,y)\mu^y_r(\dif z)\bigg|\\
&\leq C_0(1+|x|^p)e^{-\lambda (r-t)}.
\end{align*}
As a result, we obtain
\begin{align*}
|u(t,x,y)|\leq C_0(1+|x|^p)\int_t^\infty\! e^{-\lambda (r-t)}\dif r\leq C_0(1+|x|^p).
\end{align*}
Taking integral with respect to $\mu_t^y$ from both sides of (\ref{solv}) and using the property (\ref{ergo2}), we have
\begin{align*}
\int_{\mR^{d_1}}u(t,x,y)\mu_t^y(\dif x)&=\int_t^\infty\!\!\! \int_{\mR^{d_1}}P^y_{t,r}f(r,x,y)\mu_t^y(\dif x)\dif r\\
&=\int_t^\infty\!\!\!\int_{\mR^{d_1}}f(r,x,y)\mu_r^y(\dif x)\dif r=0.
\end{align*}

Next, we show that $u$ in (\ref{solv}) solves the equation (\ref{pde0}). According to \cite[Theorem 3.8]{X}, we have that for every fixed $r>0$ and $y\in\mR^{d_2}$, the function $v_r(t,x,y):=P_{t,r}f(r,x,y)\in C_p^{1+\alpha/2,2+\alpha,0}$ solves the Kolmogorov equation that
\begin{equation*}
\left\{ \begin{aligned}
&\partial_t v_r(t,x,y)+\sL_0(t,x,y)v_r(t,x,y)=0,\quad \forall(t,x)\in [0,r)\times\mR^{d_1},\\
&v_r(r,x,y)=f(r,x,y),
\end{aligned} \right.
\end{equation*}
and we have for $r-t>1$,
$$
|\p_tv_r(t,x,y)|+|\nabla_x^2v_r(t,x,y)|\leq C_0(1+|x|^p)e^{-\lambda (r-t)}.
$$
Thus, taking derivative with respect to the $t$-variable in (\ref{solv}), we obtain that
\begin{align*}
\p_tu(t,x,y)&=-f(t,x,y)+\int_t^\infty \p_tP_{t,r}f(r,x,y)\dif r\\
&=-f(t,x,y)-\int_t^\infty \sL_0(t,x,y)P_{t,r}f(r,x,y)\dif r\\
&=-f(t,x,y)-\sL_0(t,x,y)u(t,x,y).
\end{align*}
The uniqueness follows by the probability
representation  (\ref{solv}). Namely, for every solution $u$ to  (\ref{pde0}) which also satisfies the centering condition  (\ref{cen}), using It\^o's formula and taking expectation, we have for $0<s<t$ that
\begin{align*}
\mE u\big(t,X_{s,t}^y(x),y\big)&=u(s,x,y)+\mE\left(\int_s^t(\p_t+\sL_0)u\big(r,X_{s,r}^y(x),y\big)\dif r\right)\\
&=u(s,x,y)-\int_s^tP_{s,r}f(r,x,y)\dif r.
\end{align*}
Letting $t\to\infty$ and by (\ref{decay2}), we obtain (\ref{solv}).

Finally, we proceed to prove the regularity of $u$ with respect to the parameter $y$. We only prove the result for $\beta=1$,  the general case can be proved  by the same arguments as in \cite[Theorem 2.1]{RX}. We focus on the a priori estimate. Taking derivative with respect to the parameter $y$ from both sides of the equation (\ref{pde0}), we have that
\begin{align*}
&\partial_t \nabla_yu(t,x,y)+\sL_0\nabla_yu(t,x,y)\\
&\qquad\qquad\qquad=-\nabla_yf(t,x,y)-\nabla_y\sL_0u(t,x,y)=:\hat f(t,x,y),
\end{align*}
where
\begin{align*}
\nabla_y\sL_0:=\nabla_y\sL_0(t,x,y)
=\nabla_yb(t,x,y)\cdot\nabla_x+\frac{1}{2}\tr\left(\nabla_y(\sigma\sigma^*)(t,x,y)\cdot\nabla_x^2\right).
\end{align*}
The same argument as in the proof of  \cite[(28)]{P-V2} and by (\ref{FP}), we have that
$$
(\p_t+\sL_0)^*(\nabla_y\mu_t^y)=-(\nabla_y\sL_0)^*\mu_t^y,
$$
which together with the fact that $f$ satisfies (\ref{cen}) implies that
\begin{align}\label{ff}
&\int_s^t\!\!\int_{\mR^{d_1}}\hat f(r,x,y)\mu_r^y(\dif x)\dif r\no\\
&=
-\int_s^t\!\!\int_{\mR^{d_1}}\nabla_yf(r,x,y)+\nabla_y\sL_0(r,x,y)u(r,x,y)\mu_r^y(\dif x)\dif r\no\\
&=-\int_s^t\!\!\int_{\mR^{d_1}}\nabla_yf(r,x,y)\mu_r^y(\dif x)\dif r+\int_s^t\!\!\int_{\mR^{d_1}}(\p_r+\sL_0)u(r,x,y)\nabla_y\mu_r^y(\dif x)\dif r\no\\
&=-\int_s^t\nabla_y\left(\int_{\mR^{d_1}}f(r,x,y)\mu_r^y(\dif x)\right)\dif r=0.
\end{align}
Moreover, by the assumptions on the coefficients and the fact that $u(\cdot,\cdot,y)\in C_p^{1+\alpha/2,2+\alpha}$, one can check that $\hat f\in C_p^{\alpha/2,\alpha,\beta}$. Thus  we have
$$
\nabla_yu(t,x,y)=\int_t^\infty P^y_{t,r}\hat f(r,x,y)\dif r,
$$
and
$$
|\nabla_yu(t,x,y)|\leq C_0(1+|x|^p)\int_t^\infty e^{-\lambda (r-t)}\dif r\leq C_0(1+|x|^p).
$$
The proof is finished.
\end{proof}

Given a function $f(t,x,y):\mR_{+}\times\mR^{d_1}\times\mR^{d_2}\to\mR^{d_2}$, we denote by $\bar f(t,y)$ its space average with respect to $\mu_t^y$, that is,
\begin{align}\label{barf}
\bar f(t,y)=\int_{\mR^{d_1}}f(t,x,y)\mu_t^y(\dif x).
\end{align}
Concerning the regularity of the averaged function, we have the following result.

\bc\label{bf}
Assume that  $\sigma, b\in C_p^{\alpha/2,\alpha,\beta}$ with $0<\alpha\leq 1$ and $\beta\geq 0$. Given a function $f\in C_p^{\alpha/2,\alpha,\beta}$, let $\bar f$ be defined by (\ref{barf}).  Then we have $\bar f\in C_b^{\alpha/2,\beta}.$
\ec

\begin{proof}
For the regularity of $\bar f$ with respect to $y$, we only prove the result for $\beta=1$,  the general case can be proved  similarly. It is easy to see that
$$
\tilde f(t,x,y):=f(t,x,y)-\bar f(t,y)
$$
satisfies the centering condition (\ref{cen}). Using (\ref{ff}) with $f$ replaced by $\tilde f$, we have that
\begin{align*}
\int_s^t\!\!\int_{\mR^{d_1}}\!\Big[\nabla_y\tilde f(r,x,y)+\nabla_y\sL_0\tilde u(r,x,y)\Big]\mu_r^y(\dif x)\dif r=0,
\end{align*}
where $\tilde u$ is the solution of the time inhomogeneous Poisson equation (\ref{pde0}) with $f$ replaced by $\tilde f$. As a result, we obtain
\begin{align*}
\nabla_y\bar f(t,y)=\int_{\mR^{d_1}}\!\Big[\nabla_y f(t,x,y)+\nabla_y\sL_0u(t,x,y)\Big]\mu_t^y(\dif x),
\end{align*}
which implies that $\bar f(t,\cdot)\in C_b^1(\mR^{d_2})$.
As for the regularity of $\bar f$ with respect to the time variable, note that we have
\begin{align*}
\bar f(t,y)=\int_{\mR^{d_1}}P_{s,t}f(t,x,y)\mu_s^y(\dif x).
\end{align*}
Thus, $\bar f(\cdot,y)\in C_b^{\alpha/2}$ follows by the property of $P_{s,t}f(t,x,y)$. The proof is finished.
\end{proof}

\section{Weak and strong convergence in the averaging principle}

Using the technique of the infinite-horizon time-inhomogeneous Kolmogorov equation, we shall first derive two
fluctuation estimates of double averaging type for  the non-autonomous SDE (\ref{sde0})  in Subsection
4.1. Then, we prove the weak convergence in the averaging principle in Subsection 4.2. Finally, combining with the Zvonkin's transformation,   we give the proof of the strong convergence result   in Subsection
4.3.

We shall need some standard mollification arguments due to our low regularity assumptions on the coefficients. For this, let $\rho_1 : \mR\to[0,1]$ and $\rho_2 : \mR^{d_2}\to[0,1]$ be two smooth radial convolution kernel functions such that $\int_{\mR}\rho_1(r)\dif r=\int_{\mR^{d_2}}\rho_2(y)\dif y=1$, and for any $k\geq 1$, there exist constants $C_k>0$ such that $|\nabla^k\rho_1(r)|\leq C_k\rho_1(r)$ and $|\nabla^k\rho_2(y)|\leq C_k\rho_2(y)$. For every $n\in \mN^*$, set
$$
\rho^n_1(r):=n^2\rho_1(n^2 r) \qquad \text{and} \qquad \rho^n_2(y):=n^{d_2}\rho_2(ny).
$$
Given a function $f(t,s,x,y)$, define the mollifying approximations of $f$ with respect to $s$ and $y$ variables by
\begin{equation}\label{fn}
f_n(t,s,x,y):=f*\rho_2^n*\rho_1^n:=\int_0^\infty\!\!\int_{\mR^{d_2}}
f(t,s-\tilde{s},x,y-\tilde{y})\rho_2^n
(\tilde{y})\rho_1^n(\tilde{s})\dif \tilde{y}\,\dif \tilde{s}.
\end{equation}
The following  result can be proved similarly as in \cite[Lemma 4.1]{RX}, we omit the details.
\bl\label{molli}
Let $f\in C_p^{\alpha/2,\b/2,\alpha,\beta}$ with $0<\alpha,\beta\leq 2$,  and define $f_n$ by (\ref{fn}). Then we have
\begin{align}
||f(\cdot,\cdot,x,\cdot)-f_n(\cdot,\cdot,x,\cdot)||_\infty&\leq C_0\,n^{-\beta}(1+|x|^p),\label{n11}\\
||\nabla_yf_n(\cdot,\cdot,x,\cdot)||_\infty&\leq C_0\,n^{1-(\beta\wedge1)}(1+|x|^p),\label{n21}
\end{align}
and
\begin{align}\label{n22} ||\partial_tf_n(\cdot,\cdot,x,\cdot)||_\infty+||\nabla^2_yf_n
(\cdot,\cdot,x,\cdot)||_\infty\leq C_0\,n^{2-\beta}(1+|x|^p),
\end{align}
where $C_0>0$ is a constant independent of $n$.
\el

\subsection{Fluctuation estimates of double averaging type}

Recall that $\{\mu_t^y\}_{t\geq 0}$ is  an evolution system of invariant measures for the frozen process (\ref{froz}). Given a function $f:\mR_{+}^2\times\mR^{d_1+d_2}\to\mR^{d_2}$,
we say  that $f$ satisfies  the assumption $\bf{(H_f)}$ if the following condition holds:

\vspace{2mm}
\begin{enumerate}[$ {\bf ({H_f})} $]
\item there exists a function $\bar {\bar f}:\mR_{+}\times\mR^{d_2}\to\mR^{d_2}$ such that for all $T>0$ and $y\in\mR^{d_2}$,
\begin{align}\label{bar-f4}
\left|\frac{1}{T}\int_{0}^{T}\!\!\int_{\mR^{d_1}}f(t,s,x,y)\mu_t^y(\dif x)\dif t-\bar{\bar f}(s,y)\right|\leq \kappa(T),
\end{align}
where $\kappa:\mR_+\to\mR_+$ satisfies $\kappa(T)\to 0$ as $T\to \infty.$
\end{enumerate}

For simplicity, we denote by
\begin{equation}\label{pu}
\partial_1 f(t,s,x,y):=\partial_t f(t,s,x,y),\quad \partial_2 f(t,s,x,y):=\partial_s f(t,s,x,y),
\end{equation}
and
\begin{align}\label{L1}
\cL_y:=\cL_y(t/\eps,x,y)=F(t/\eps,x,y)\cdot\nabla_y+\frac{1}{2}
\tr\left(GG^*(t/\eps,x,y)\cdot\nabla_y^2\right).
\end{align}
The following result will be critical  for proving the weak convergence in the averaging principle.
\bl[Weak fluctuation estimate]\label{xx}
Let $T>0$ and $\bf{(H_1)}$ hold. Assume that $\sigma, b\in C_p^{\alpha/2,\alpha,\beta}$   with $\alpha,\beta>0$ and $F,G\in L_p^\infty$. Then for every $f\in C_p^{\alpha/2,\beta/2,\alpha,\beta}$ satisfying $\bf{(H_f)}$ and $t\in [0,T],$
we have
\begin{align*}
\left|\mE\bigg(\int_0^t \big[f(r/{\eps},r,X_r^{\eps},Y_r^{\eps})-\bar{\bar f}(r,Y_r^\eps)\big]\dif r\bigg)\right|\leq C_T\,\big(\eps^{(\beta/2)\wedge1}+t\cdot \kappa(t/{\eps})\big),
\end{align*}
where $C_T>0$ is a constant, and $\bar{\bar{f}}$ is given in (\ref{bar-f4}).
\el
\begin{proof}
We write
\begin{align*}
&\left|\mE\bigg(\int_0^t\big[f\big(r/\eps,r
,X_r^\eps,Y_r^\eps\big)
-\bar{\bar{f}}(r,Y_r^\eps)\big]\dif r\bigg)\right|\no\\
&\leq\left|\mE\bigg(\int_0^t\big[f\big(r/\eps,r,X_r^\eps,Y_r^
\eps\big)-\bar{f}(r/\eps,r,Y_r^\eps)\big]\dif r\bigg)\right
|\no\\
&+\left|\mE\bigg(\int_0^t\big[\bar{f}\big(r/\eps,r,Y_r^\eps\big)-
\bar{\bar{f}}(r,Y_r^\eps)\big]\dif r\bigg)\right|
=:\mathcal T_1(t,\eps)+\mathcal T_2(t,\eps),
\end{align*}
where $\bar f(t,s,y)$ is defined by
\begin{align*}
\bar f(t,s,y)=\int_{\mR^{d_1}}f(t,s,x,y)\mu_t^y(\dif x).
\end{align*}
Without loss of generality, we may assume that $\b\leq 2$. Below, we proceed to control each term separately.

\vspace{1mm}
\noindent {\it (i) Estimate of $\cT_1(t,\eps)$}. Consider the following time-inhomogeneous Kolmogorov equation in $[0,\infty)\times\mR^{d_1}:$
\begin{align}\label{pde1}
&\partial_t \psi(t,s,x,y)+\cL_0\psi(t,s,x,y)=-(f(t,s,x,y)-\bar{f}(t,s,y)),
\end{align}
where $\sL_0$ is given by (\ref{L0}), and $(s,y)\in\mR_+\times\mR^{d_2}$ are regarded as parameters.
Note that $f(t,s,x,y)-\bar{f}(t,s,y)$ satisfies the centering condition (\ref{cen}):
$$
\int_{\mR^{d_1}}(f(t,s,x,y)-\bar{f}(t,s,y))\mu_t^y(\dif x)=0,\quad\forall (s,y)\in\mR_+\times\mR^{d_2}.
$$
Moreover,  under the assumptions on the coefficients and by Corollary \ref{bf}, we have  $\bar{f}(t,s,y)\in C_b^{\a/2,\b/2,\b}$. By Theorem \ref{cauchy1} there exists a unique solution $\psi\in C_p^{1+\alpha/2,\beta/2,2+\alpha,\beta}$  to equation (\ref{pde1}).
Let $\psi_n$ be the mollifier of $\psi$ defined by (\ref{fn}), i.e.,
\begin{align*}
\psi_n(t,s,x,y):=\psi*\rho_2^n*\rho_1^n:=\int_0^\infty\!\!\int_{\mR^{d_2}}
\psi(t,s-\tilde{s},x,y-\tilde{y})\rho_2^{n}(\tilde y)\rho_1^{n}(\tilde s)\dif \tilde y\dif \tilde s.
\end{align*}
Applying It\^o's formula to  $\psi_n\big(t/{\eps},t,X_t^\eps,Y_t^\eps\big)$,
 we have
\begin{align*}
\psi_n(t/\eps,t,X_t^\eps,Y_t^\eps)&=\psi_n(0,0,x,y)+\frac{1}
{\eps}\left(\int_0^t(\partial_1+\cL_0)\psi_n(r/\eps,r,X_r^\eps,
Y_r^\eps)\dif r\right)\no\\
&+\int_0^t(\p_2+\cL_y)\psi_n(r/\eps,r,X_r^\eps,Y_r^\eps)\dif r\no\\
&+\int_0^t\nabla_y\psi_n(r/\eps,r,X_r^\eps,Y_r^\eps)G(r/\eps,
X_r^\eps,Y_r^\eps)\dif W_r^2\no\\
&+\frac{1}{\sqrt{\eps}}\int_0^t\nabla_x\psi_n(r/\eps,r,X_r^
\eps,Y_r^\eps)\sigma(r/\eps,X_r^\eps,Y_r^\eps)\dif W_r^1 ,
\end{align*}
where the operator $\sL_y$ is defined by (\ref{L1}), and $\p_1\psi_n,\p_2\psi_n$ are given by (\ref{pu}). Multiplying  both sides of the above equality by $\eps$ and combining with (\ref{pde1}), we obtain
\begin{align}\label{f-barf}
&\int_0^t[f(r/\eps,r,X_r^\eps,Y_r^\eps)-\bar f(r/\eps,r,
Y_r^\eps)]\dif r\no\\
&=\eps[\psi_n(0,0,x,y)-\psi_n(t/\eps,t,X_t^\eps,Y_t^\eps)]
+\eps\int_0^t(\p_2+\sL_y)\psi_n(r/\eps,r,X_r^\eps,Y_r^\eps)\dif r\no\\
&\quad+\eps\int_0^t\nabla_y\psi_n(r/\eps,r,X_r^\eps,Y_r^\eps)
G(r/\eps,X_r^\eps,Y_r^\eps)\dif W_r^2\no\\
&\quad+\sqrt{\eps}\int_0^t\nabla_x\psi_n(r/\eps,r,X_r^\eps,
Y_r^\eps)\sigma(r/\eps,X_r^\eps,Y_r^\eps)\dif W_r^1\no\\
&\quad+\int_0^t(\partial_1+\sL_0)(\psi_n-\psi)(r/\eps,r,X_r
^\eps,Y_r^\eps)\dif r.
\end{align}
Taking expectation, we get
\begin{align}\label{Zi}
&\mathcal T_1(t,\eps)
\leq\eps\mE\left|\psi_n(0,0,x,y)-\psi_n\big(t/{\eps}
,t,X_t^
\eps,Y_t^\eps\big)\right|\no\\
&\qquad\qquad+\eps\mE\left|\int_0^t(\partial_2+\cL_y) \psi_n\big(r/{\eps},r,
X_r^\eps,Y_r^\eps\big)\dif r\right|\no\\
&\qquad\qquad+\mE\left|\int_0^t(\partial_1+\cL_0)(\psi_n-\psi)
\big(r/{\eps},r,X_r^\eps,Y_r^\eps\big)\dif r\right|=:\sum_{i=1}^3\cT_{1,i}(t,\eps).
\end{align}
Under the assumptions on the coefficients, we have (see e.g. \cite[Proposition 1]{P-V} and \cite[Lemma 3.1]{G}) that for any $q > 0$,
\begin{align}\label{ex}
\sup\limits_{\eps\in(0,1)}\sup_{t\in[0,T]}\mE\big(|X_t^\eps |^q+|Y_t^\eps |^q\big) \leq C_{q,T}<\infty.
\end{align}
Since $\psi\in C_p^{1+\alpha/2,
\beta/2,2+\alpha,\beta}$ , we deduce that there exists a $p>1$ such that for any $t\in[0,T]$,
$$
\cT_{1,1}(t,\eps)\leq C_{T}\,\eps(1+\mE|X_t^
\eps|^p)\leq C_{T}\,\eps.
$$
By (\ref{n21}) and (\ref{n22}) of Lemma \ref{molli} and the assumptions that $F,G\in L_p^\infty(\mR_+\times\mR^{d_1+d_2})$, we have that there exists $p>1$ such that
\begin{align*}
&|(\partial_2+\cL_y ) \psi_n(t,s,x,y)|\\
&=\Big|\Big(\partial_2\psi_n+\frac{1}{2}\tr\big
(GG^\ast(t,x,y)\nabla_y^2\psi_n\big)+F(t,x,y)
\!\cdot\!\nabla_y\psi_n\Big)(t,s,x,y)\Big|\\
&\leq C_{T}\big(1+|x|^p\big)\Big(|\nabla_y
\psi_n(t,s,x,y)|+ |\nabla_y^2 \psi_n(t,s,x,y)|+|\partial_2\psi_n(t,s,x,y)|\Big)\\
&\leq C_{T}n^{2-\beta}\big(1+|x|^{2p}\big).
\end{align*}
As a result, we have for any $t\in[0,T]$,
\begin{align*}
\cT_{1,2}(t,\eps)&\leq C_{T}\,\eps n^{2-\beta}\int_0^T\Big(1+\mE|X_r^\eps|^
{2p}\Big)\dif r\\
&\leq C_{T}\,\eps\, n^{2-\beta}.
\end{align*}
For the last term, by the fact that
$$
\partial_1(\psi_n)=(\partial_1 \psi)\ast\rho_2^n\ast\rho_1^n,\quad
\nabla_x^2(\psi_n)=(\nabla_x^2 \psi)\ast \rho_2^n\ast\rho_1^n,
$$
and using (\ref{n11}), we derive that for any $t\in[0,T]$,
\begin{align*}
\cT_{1,3}(t,\eps)&\leq C_{T}\mE\Big(\int_0^T|(
\partial_1 \psi_n-\partial_1 \psi)\big(r/{\eps},
r,X_r^\eps,Y_r^\eps\big)|\dif r\Big)\\
&~~~~+C_{T}\mE\Big(\int_0^T\big(1+|X_r^\eps|^p
\big)\sum_{i=1,2}\Big|(\nabla_x^i \psi_n
-\nabla_x^i \psi)\big(r/{\eps},r,X_r^\eps,Y_r^\eps\big)\Big|\dif r\Big)\\
&\leq C_{T}\,n^{-\beta}\int_0^T\big(1+\mE|X_r^
\eps|^{2p}\big)\dif r\leq C_{T}\,n^{-\beta}.
\end{align*}
Taking $n=\eps^{-1/2}$, we obtain for any $t\in[0,T]$,
$$
\cT_1(t,\eps)\leq C_{T}\big(\eps+\eps n^{2-\beta}+n^{-\beta}\big)\leq C_{T}\,\eps^{\beta/2}.
$$

\vspace{1mm}
\noindent {\it (ii) Estimate of $\cT_2(t,\eps)$}.
Let us consider the following ordinary differential equation:
\begin{equation}\label{ode}
\left\{\begin{aligned}
&\partial_t U(t,s,y)=-\big(\bar{f}(t,s,y)-\bar{\bar{f}}(s,y)\big),\\
&U(0,0,y)=0,
\end{aligned}\right.
\end{equation}
where $(s,y)\in \mR_+\times\mR^{d_2}$ are regarded as parameters. The solution $U(t,s,y)$ is given by
\begin{equation}\label{U}
U(t,s,y)=-\int_0^t\big(\bar{f}(r,s,y)-\bar{\bar{f}}
(s,y)\big)\dif r.
\end{equation}
Let $ U_n(t,s,y)$ be the mollifier of $ U(t,s,y)$ as in (\ref{fn}), i.e.,
\begin{align*}
U_n(t,s,y):=U*\rho_2^n*\rho_1^n:=\int_0^\infty\!\!\int_{\mR^{d_2}}
U(t,s-\tilde{s},y-\tilde{y})\rho_2^{n}(\tilde y)\rho_1^{n}(\tilde s)\dif \tilde y\dif \tilde s.
\end{align*}
Applying It\^o's formula to  $U_n\big(t/{\eps},t,Y_t^\eps\big)$, we have
\begin{align}\label{Un}
U_n\big(t/{\eps},t,Y_t^\eps\big)
&=\frac{1}{\eps}\int_0^t\partial_1U_n\big(r/{\eps}
,r,Y_r^\eps\big)\dif r+\int_0^t(\partial_2+\cL_y)U_n\big(r/{\eps},
r,Y_r^\eps\big)\dif r\no\\
&+\int_0^tG(r/\eps,X_r^\eps,Y_r^\eps)
\nabla_yU_n\big(r/{\eps},r,Y_r^\eps\big)\dif W_r^2,
\end{align}
where $\p_1U_n, \p_2U_n$ and $\sL_yU_n$ are given by (\ref{pu}) and (\ref{L1}), respectively.
Multiplying  both sides of the above equality by $\eps$, taking expectation by (\ref{ode}) and using arguments similar to those in (\ref{Zi}), we obtain for any $t\in[0,T]$,
\begin{align}\label{Unn}
&\left|\mE\int_0^t\big[\bar{f}\big(r/
{\eps},r,Y_r^\eps\big)-\bar{\bar{f}}(r,Y_r^\eps)
\big]\dif r\right|\no\\
&\leq\eps\mE\Big|U_n\big(t/
{\eps},t,Y_t^\eps\big)\Big|+\eps
\mE\Big|\int_0^t(\partial_2+\cL_y)
U_n\big(r/{\eps},r,Y_r^\eps\big)
\dif r\Big|\no\\
&+\mE\Big|\int_0^t\partial_1(U_n-U)
\big(r/{\eps},r,Y_r^\eps\big)\dif r
\Big|\no\\
&\leq\eps\mE\Big|U_n\big(t/{\eps},t,Y_t^\eps\big
)\Big|+C_{T}
\big(\eps n^{2-\b}+n^{-\b}\big).
\end{align}
By (\ref{U}) and the assumption (\ref{bar-f4}), we deduce that
\begin{equation*}
\big|U\big(t/{\eps},t,y\big)\big|\leq\frac{t}{\eps}
\bigg|\frac{1}
{t/\eps}\int_0^{t/\eps}\big[\bar{f}(r,s,y)-
\bar{\bar{f}}(s,y)\big]\dif r\bigg|\leq \eps^{-1}\,t\cdot \kappa\big(t/{\eps}\big).
\end{equation*}
Thus we have for $t\in [0,T]$,
$$
\eps\mE\big|U_n\big(t/{\eps},t,Y_t^\eps\big)
\big|\leq \eps\mE\Big|U\big(t/{\eps},t,Y_t^\eps\big)\Big|
\leq C_{T}~t\cdot \kappa\big(t/{\eps}\big).
$$
Taking $n=\eps^{-1/2}$, we obtain
$$
\cT_2(t,\eps)= C_{T}\Big(t\cdot\kappa\big(t/{\eps}
\big)+\eps n^{2-\b}+n^{-\b}\Big)\leq C_{T}\Big(\eps^{\b/2}+t\cdot \kappa\big(t/{\eps}\big)\Big).
$$
Combining the above computations, we obtain the desired result.
\end{proof}

Next, we proceed to derive the  strong fluctuation estimate for the non-autonomous multi-scale systems (\ref{sde0}), which will be used to prove the strong convergence in the averaging principle. Unlike Lemma \ref{xx}, we shall only need to handle function $f:\mR_{+}\times\mR^{d_1+d_2}\to\mR^{d_2}$ which does not  depend on the $s$-variable. Thus, we say a function $f(t,x,y)$ satisfies $\bf{({H}_f)}$ if there exists a function $\bar{\bar f}(y)$ which is independent of $s$-variable such that (\ref{bar-f4}) holds.

\bl[Strong fluctuation estimate]\label{xxx}
Let $T>0$ and $\bf{({H}_1)}$ hold. Assume that $\sigma,b\in C_p^{\alpha/2,\alpha,\beta}$ with $\alpha,\beta>0$, $F,G\in L_p^\infty$, and $f\in C_p^{\alpha/2,\alpha,\beta}$ satisfies $\bf{({H}_f)}$. Then for any $t\in[0,T]$ and $q\geq 2,$
\begin{align*}
\mE\left|\int_0^t\big[f(r/\eps,X_r^\eps,Y_r^\eps)-\bar {\bar{f}}(Y_r^\eps)\big] \dif r\right|^q\leq C_T\Big(\eps^{(\beta\wedge1) q/2}+\big[t\cdot\kappa(t/\eps)\big]^q\Big),
\end{align*}
where $C_T>0$ is a constant independent of $\beta$ and $\eps$.
\el
\begin{proof}
For every $q\geq 2,$ we write that
\begin{align*}
&\mE\left|\int_0^t\big[f(r/\eps,X_r^\eps,Y_r^\eps)-\bar {\bar{f}}(Y_r^\eps)\big]\dif r\right|^q\no\\&\leq C_q\mE\bigg|\int_0^t\big[f(r/\eps,X_r^\eps,Y_r^\eps)-  \bar f(r/\eps,Y_r^\eps)\big]
\dif r\bigg|^q\no\\
&+ C_q\mE\bigg|\int_0^t\big[\bar f(r/\eps,Y_r^\eps)-\bar{\bar f}(Y_r^\eps)\big]\dif r\bigg|^q=:\sS_1(t,\eps)+\sS_{2}(t,\eps),
\end{align*}
where $\bar f(t,y)$ is defined  by (\ref{barf}).
Below, we assume that $\b\leq 1$. To control the above two terms, we give the main difference between the proof of Lemma \ref{xx}.

\vspace{1mm}
\noindent {\it (i) Estimate of $\sS_1(t,\eps)$}. Let $\psi$ be the unique solution of (\ref{pde1}). Note that $\psi$ does not depend on $s$-variable, i.e.,
\begin{equation*}
\partial_t\psi(t,x,y)+\sL_0\psi(t,x,y)=-[f(t,x,y)-  \bar f(t,y)],
\end{equation*}
where $\sL_0$ is given by (\ref{L0}), and $y\in \mR^{d_2}$ is regarded as  a parameter.
Following the same arguments as in (\ref{f-barf}) and (\ref{Zi}), we have
\begin{align*}
\sS_1(t,\eps)&\leq C_q{\eps^q}\big[\left|\psi_n(0,x,y)\right|^q+\mE
\left|\psi_n(t/\eps,X_t^\eps,Y_t^\eps)\right|^q\big]\\
&\quad +C_q\eps^q\mE\left|\int_0^t\sL_y\psi_n(r/\eps,X_r^\eps,Y_r^\eps)\dif r\right|^q\\
&\quad +C_q\mE\left|\int_0^t(\partial_r+\sL_0)(\psi_n-\psi)(r/\eps,X_r^\eps,Y_r^\eps)\dif r\right|^q\\
&\quad+C_q\eps^q\mE\left|\int_0^t\nabla_y\psi_n(r/\eps,X_r^\eps,Y_r^\eps)
G(r/\eps,Y_r^\eps)\dif W_r^2\right|^q\\
&\quad+C_q\eps^{q/2}\mE\left|\int_0^t\nabla_x\psi_n(r/\eps,X_r^\eps,Y_r^\eps)
\sigma(r/\eps,X_r^\eps,Y_r^\eps)\dif W_r^1\right|^q=:\sum_{i=1}^5\sS_{1,i}(t,\eps).
\end{align*}
The first three terms can be estimated exactly as in the estimates for (\ref{Zi}), and we have
$$
\sS_{1,1}(t,\eps)+\sS_{1,2}(t,\eps)+\sS_{1,3}(t,\eps)\leq C_T(\eps+\eps n^{2-\beta}+n^{-\b})^q.
$$
We only need to deal with the martingale terms $\sS_{1,4}$ and $\sS_{1,5}$. By (\ref{n22}),  (\ref{ex}) and the H\"older inequality, we have
\begin{align*}
\sS_{1,4}(t,\eps)&\leq C_{T}\,\eps^qn^{q(1-\beta)} \mE\left(\int_0^t\big(1+|X_r^\eps|^{2p}\big)\dif r\right)^{q/2}\\
&\leq C_{T}\,\eps^qn^{q(1-\beta)} \mE\int_0^t\left(1+|X_r^\eps|^{pq}\right)\dif r\leq C_{T}\,\eps^qn^{q(1-\beta)}, \end{align*}
and similarly,
$$\sS_{1,5}(t,\eps)\leq C_{T}\,\eps^{q/2} \mE\int_0^t\left(1+|X_r^\eps|^{pq}\right)\dif r\leq C_{T}\,\eps^{q/2 }.$$
Taking $n=\eps^{-1/2},$ we obtain for any $t\in[0,T]$,
$$\sS_{1}(t,\eps)\leq C_{T}\,\big(\eps^{q/2}+\eps^qn^{q(2-\beta)}+n^{-\beta q}\big)\leq C_{T}\,\eps^{\beta q/2}. $$

\vspace{1mm}
\noindent {\it (ii) Estimate of $\sS_2(t,\eps)$}.
Let $\bar U(t,y)$ solve (\ref{ode}) with the right hand side replaced by $\bar f(t,y)-{\bar{\bar f}}(y)$. Similar to (\ref{Un}) and (\ref{Unn}),  we have for  any $t\in[0,T]$ and $q\geq 2,$
\begin{align}\label{UI}
&\mE\bigg|\int_0^t [\bar f(r/\eps,Y_r^\eps)-\bar{\bar f}(Y_r^\eps)]\dif r\bigg|^q\no\\
&\leq C_q\bigg[{\eps^q}\bigg(\mE|\bar U_n(t/\eps,Y_t^\eps)|^q+\mE\left|\int_0^t
\sL_y
\bar U_n(r/\eps,Y_r^\eps)\dif r\right|^q\bigg)\no\\
&+\mE\left|\int_0^t\partial_r(\bar U_n-\bar U)(r/\eps,Y_r^\eps)\dif r\right|^q+\eps^q\mE\left|\int_0^t\nabla_y\bar U_n(r/\eps,Y_r^\eps)G(r/\eps,Y_r^\eps)\dif W_r^2\right|^q\bigg]\no\\
&
\leq C_{q,T}\Big(\eps n^{(2-\beta)}
+n^{-\beta}+t\cdot\kappa(t/\eps)\Big)^q\no\\
&+C_q\,{\eps^q} \mE\left|\int_0^t\nabla_y\bar U_n(r/\eps,Y_r^\eps)G(r/\eps,Y_r^\eps)\dif W_r^2\right|^q.
\end{align}
For the martingale term, we have by (\ref{n21}) that
$$
\mE\left|\int_0^t\nabla_y\bar U_n(r/\eps,Y_r^\eps)G(r/\eps,Y_r^\eps)\dif W_r^2\right|^q\leq C_{q,T}n^{(1-\beta)q}.
$$
Taking $n=\eps^{-1/2},$ we obtain
\begin{align*}
\sS_2(t,\eps)
\leq C_{T}\Big(\eps^{\beta q/2}+\big[t\cdot\kappa(t/\eps)\big]^q\Big). \end{align*}
Combining the above computations, we obtain the desired result.
\end{proof}

\subsection{Proof of Theorem \ref{main1}}
Recall that $\bar{Y}_t$ satisfies the SDE (\ref{wbsde}). We shall also write $\bar Y_t(y)$  to stress its dependence on the
initial value.
Let $\bar{\cL}$ be the  generator of $\bar{Y}_t$, i.e.,
\begin{equation}\label{barL}
\bar{\cL}:=\bar{\cL}(y):=\bar{\bar F}(y)\cdot\nabla_y+\frac{1}{2}\tr\big(\bar{\bar G}\bar{\bar G}^*(y)\cdot\nabla_y^2\big),
\end{equation}
where $\bar{\bar F}$ and $\bar{\bar G}$ are given by (\ref{bF}) and (\ref{bGG}), respectively. Consider the following  Kolmogorov equation in $\mR_+\times \mR^{d_2}$:
\begin{equation}\label{PDE-B}
\left\{\begin{aligned}
&\partial_{t}u(t,y) -\bar{\cL}u(t,y)=0,\\
&\ u(0,y) =\varphi(y),
\end{aligned} \right.
\end{equation}
where $\varphi$ is a measurable function. We have the following result.
\bl
Suppose that the assumptions in Theorem \ref{main1} hold. Then for every $\varphi\in C_b^{2+\b}$, there exists a unique solution $u\in C_b^{(2+\beta)/2,2+\beta}$ to equation (\ref{PDE-B}) which is given by
\begin{align}\label{cau0}
u(t,y)=\mE\varphi\big( \bar{Y}_t(y)\big).
\end{align}
\el
\begin{proof}
Under our assumptions and by Corollary \ref{bf}, one can check that
$
\bar{\bar F}, \bar{\bar G}\in C_b^\beta(\mR^{d_2})$
  and $\bar{\bar G}$ is   uniformly elliptic. By \cite[Theorem 3.4]{RX1}, we have the desired result.
\end{proof}
Now, we are in the position to give:

\begin{proof}[Proof of Theorem \ref{main1}]
Fix $T>0$ and $\varphi\in C_b^{2+\b}$, let $u$ be given by (\ref{cau0}). Define
$$
\tilde{u}(t,y):=u(T-t,y).
$$
Then it is obvious that
$$
\tilde{u}(0,y)=\mE\varphi(\bar{Y}_T(y)),\quad
\tilde{u}(T,y)=\varphi(y),
$$
and by (\ref{PDE-B})
\begin{equation}\label{PDE-L}
\partial_{t}\tilde{u}(t,y) =-\bar{\cL}\tilde{u}(t,y).
\end{equation}
As a result, we can deduce by (\ref{barL}) (\ref{cau0}), (\ref{PDE-L}) and It\^o's formula that
\begin{align}\label{cI}
&\cI(T,\eps):=\big|\mE\varphi(Y_T^{\eps})-\mE
\varphi(\bar{Y}_T)\big|=\big|\mE\tilde{u}(T,Y_T^
\eps)-\mE \tilde{u}(0,y)\big|\no\\
&=\bigg|\mE\bigg(\int_0^T\partial_r\tilde{u}(r,
Y_r^\eps)+\frac{1}{2}\tr\Big
(GG^\ast\big(r/{\eps},X_r^\eps,Y_r^\eps\big)\cdot
\nabla_y^2 \tilde{u}(r,Y_r^
\eps)\Big)\no\\
&\qquad\qquad\qquad\qquad\quad\,+F\big(r/{\eps},X_r^\eps,Y_r^\eps\big)\cdot
\nabla_y \tilde{u}(r,Y_r^\eps)\dif r\bigg)\bigg|\no\\
&\leq\left|\mathbb{E}\bigg(\int_0^T\frac{1}{2}\tr\Big(\big[
GG^*(r/{\eps},X_r^\eps,Y_r^\eps)-\bar{\bar G}\bar{\bar G}^*(Y_r^\eps)\big]\cdot\nabla_y^2\tilde u(r,Y_r^\eps)\Big)\mathrm{d}r\bigg)\right|\no\\
&+\left|\mathbb{E}\Big(\!\!\int_0^T\Big[F(r/{\eps},X_r^\eps,
Y_r^\eps)-\bar{\bar F}(Y_r^\eps)\Big]\cdot\nabla_y\tilde u(r,Y_r^\eps)\mathrm{d}r\Big)\right|=:\cI_1(T,\eps)+
\cI_2(T,\eps).
\end{align}
Let
$$
f(t,s,x,y):=\frac{1}{2}\tr\big(
GG^*(t,x,y)\cdot\nabla_y^2\tilde u(s,y)\big).
$$
Then, by the assumption (\ref{bGG}) one can check that $f$ satisfies $(\mathbf{ H_f})$ with
$$
\bar{\bar{f}}(s,y)=\frac{1}{2}\tr\big(\bar{\bar{G}}\bar{\bar{G}}^*(y)
\cdot\nabla_y^2\tilde u(s,y)\big).
$$
Consequently, as a direct result of the Lemma \ref{xx}, we get
$$
\cI_1(T,\eps)\leq C_{T}\Big(\eps^{(\b/2)\wedge 1}+T\cdot \kappa_2\big(T/{\eps}\big)\Big).
$$
Similarly, by the assumption (\ref{bF}) and Lemma \ref{xx}, we have
$$
\cI_2(T,\eps)\leq C_{T}\Big(\eps^{(\b/2)\wedge 1}+T\cdot \kappa_1\big(T/{\eps}\big)\Big).
$$
Thus we obtain
$$
\cI(T,\eps)\leq C_{T}\Big(\eps^{(\b/2)\wedge 1}+T\cdot \kappa_1\big(T/{\eps}\big)+T\cdot \kappa_2\big(T/{\eps}\big)\Big).
$$
The proof is finished.
\end{proof}

\subsection{Proof of Theorem \ref{main11}}

It seems to be difficult to prove the strong convergence of $Y_t^\eps$ to $\bar{Y}_t$ directly due to the low regularity (only H\"older continuous) of the drift coefficient $F$. For this reason, we shall use the Zvonkin's transformation as in  \cite{RX1} to transform the equations of $Y_t^\eps$ and $\bar Y_t$ into new ones.

Below, we assume the assumptions in Theorem \ref{main11} hold. For $\lambda>0$, consider the following elliptic equation in $\mR^{d_2}$:
\begin{align} \label{ellip}
\lambda u(y)- \bar\sL_y u(y)=\bar{\bar F}(y),
\end{align}
where $ \bar\sL_y$ is defined by (\ref{barL}) with $\bar{\bar G}$ replaced by $\bar G$, i.e.,
\begin{align}\label{bar1}
\bar\sL_y:=\bar\sL(y):=\bar{\bar F}(y)\cdot\nabla_y+\frac{1}{2}\tr\left({ \bar G}{ \bar G}^*(y)\cdot\nabla^{2}_y\right).
\end{align}
Since $\bar{\bar F}\in C_b^\b(\mR^{d_2})$
and $\bar G\in C_b^1(\mR^{d_2})$  is uniformly elliptic,  there exists a unique solution $u\in C_b^{2+\beta}(\mR^{d_2})$ to equation (\ref{ellip})  (see e.g.  \cite[Chapter IV, Section 5]{LSU}). Moreover, we can take $\lambda$ large enough such that
\begin{align}\label{derit}
|\nabla_y u(y)|\leq 1/2, \quad\forall \,y\in \mR^{d_2}.
\end{align}
Define two new processes by
\begin{align}\label{byy}  \bar V_t:=  \bar Y_t+u(\bar Y_t)\end{align}
and
\begin{align} \label{yy} V_t^\eps:=  Y_t^\eps+u(Y_t^\eps).
\end{align}
We have the following results.

\bl[Zvonkin's transformation]
Let $  \bar V_t$ and $V_t^\eps$  be defined by (\ref{byy}) and (\ref{yy}), respectively. Then we have
\begin{align}\label{bary}
\dif   \bar V_t= \lambda u(  \bar Y_t)\dif t+\bar G(  \bar Y_t)(\mI+\nabla_y u)(\bar Y_t)\dif W_t^2, \, \, \bar V_0=y+u(y),
\end{align}
and
\begin{align}\label{xux}
\dif V_t^\eps=&\lambda u(Y_t^\eps)\dif t+G(t/\eps,Y_t^\eps)(\mI+\nabla_yu)(Y_t^\eps)\dif W_t^2\no\\
&+\big[F(t/\eps,X_t^\eps,Y_t^\eps)-\bar{\bar F}(Y_t^\eps)\big](\mI+\nabla_yu)(Y_t^\eps)\dif t\no\\
&+\frac{1}{2}\tr\Big(\big(GG^*(t/\eps,Y_t^\eps)-\bar G\bar G^*(Y_t^\eps)\big)\cdot\nabla^2_y u(Y_t^\eps)\Big)\dif t, \,\, V_0^\eps=y+u(y).
\end{align}
\el

\begin{proof}
We only prove the formula (\ref{xux}) since the proof of (\ref{bary}) is easier and follows by the same argument.
Applying It\^o's formula to $u(Y_t^\eps),$ we have
\begin{align*}
u(Y_t^\eps)&=u(y)+\int_0^t\sL_yu(Y_r^\eps)\dif r+\int_0^tG(r/\eps,Y_r^\eps)\nabla_yu(Y_r^\eps)\dif W_r^2\no\\
&=u(y)+\int_0^t\bar\sL_y u(Y_r^\eps)\dif r+\int_0^tG(r/\eps,Y_r^\eps)\nabla_yu(Y_r^\eps)\dif W_r^2\no\\
&\quad+\int_0^t\big[F(r/\eps,X_r^\eps,Y_r^\eps)-\bar{\bar F}(Y_r^\eps)\big]\cdot\nabla_y u(Y_r^\eps)\dif r\\
&\quad+\frac{1}{2}\int_0^t \tr\Big(\big(GG^*(r/\eps,Y_r^\eps)-\bar G\bar G^*(Y_r^\eps)\big)\cdot\nabla^2_y u(Y_r^\eps)\Big)\dif r\no\\
&=u(y)+\int_0^t\lambda u(Y_r^\eps)-\bar{\bar F}(Y_r^\eps)\dif r+\int_0^tG(r/\eps,Y_r^\eps)\nabla_yu(Y_r^\eps)\dif W_r^2\no\\&\quad+\int_0^t\big[F(r/\eps,X_r^\eps,Y_r^\eps)-\bar{\bar F}(Y_r^\eps)\big]\cdot\nabla_y u(Y_r^\eps)\dif r\no\\&\quad+\frac{1}{2}\int_0^t\tr\Big(\big(GG^*(r/\eps,Y_r^\eps)-\bar G\bar G^*(Y_r^\eps)\big)\cdot\nabla^2_y u(Y_r^\eps)\Big)\dif r.
\end{align*}
Adding this with the equation of $Y_t^\eps$,
we obtain the desired result.
\end{proof}

Now, we are in the position to give:
\begin{proof}[Proof of Theorem \ref{main11}]
For simplicity, we prove the result for $q=2$. Recall that $u$ solves the equation (\ref{ellip}).
Taking $\lambda>0$   such that (\ref{derit}) holds, then for every $t\in[0,T]$,
\begin{align*}
\mE|Y_{t}^\eps-  \bar Y_{t}|^2\leq C_0\,\mE| V^\eps_{t}-\bar V_{t}|^2.
\end{align*}
In view of  (\ref{bary}) and (\ref{xux}), we deduce that
\begin{align}\label{mar}
&\mE\big|Y_{t}^\eps-  \bar Y_{t} \big|^2\leq
C_0\mE\bigg|\int_0^{t}\lambda\left(u(Y_r^\eps)-u(\bar Y_r)\right)\dif r\bigg|^2\no\\
&+
C_0\mE\bigg|\int_0^{t}\big[\bar G(Y_r^\eps)(\mI+\nabla_yu)(Y_r^\eps)-
\bar G( \bar Y_r)(\mI+\nabla_yu)(\bar Y_r)\big]\dif W_r^2\bigg|^2\no\\
&+C_0\mE\bigg|\int_0^{t}\big[G(r/\eps,Y_r^\eps)-
\bar G(Y_r^\eps)\big](\mI+\nabla_y u)(Y_r^\eps)\dif W_r^2\bigg|^2\no\\
&+C_0\mE\bigg|\int_0^{t}[F(r/\eps,X_r^\eps,Y_r^\eps)-  \bar{\bar F}(Y_r^\eps)](\mI+\nabla_yu)(Y_r^\eps)
\dif r\bigg|^2\no\\
&+C_0\mE\bigg|\int_0^{t}\tr\Big(\big(GG^*(r/\eps,Y_r^\eps)-\bar G\bar G^*(Y_r^\eps)\big)\cdot\nabla^2_y u(Y_r^\eps)\Big)\dif r\bigg|^2.
\end{align}
By the fact that $u\in C_b^{2+\beta}(\mR^{d_2})$, $\bar G\in C_b^1(\mR^{d_2})$ and the estimate (\ref{ex}), we further obtain that for every $t\in[0,T]$,
\begin{align*}
\mE\big|Y_{t}^\eps-  \bar Y_{t} \big|^2&\leq C_{\lambda,T}\int_0^t\mE|Y_r^\eps-\bar Y_r|^2\dif r+C_{T}\,\mE\left(\int_0^t|G(r/\eps,Y_r^\eps)-
\bar G( Y_r^\eps)|_{HS}^2\dif r\right)\\
&\quad+C_{\lambda,T}\,\mE\bigg|\int_0^{t}\big[F(r/\eps,X_r^\eps,Y_r^\eps)-  \bar{\bar F}(Y_r^\eps)\big]\cdot(\mI+\nabla_yu)(Y_r^\eps)
\dif r\bigg|^2.
\end{align*}
Let
$$
f(t,x,y):=F(t,x,y)(\mI+\nabla_y u(y))\quad\text{and}\quad \bar{\bar f}(y):=\bar{\bar F}(y)(\mI+\nabla_y u(y)).
$$
By the assumption (\ref{bF}), one can check that $\bf{( H_f)}$ holds.  Applying Lemma \ref{xxx} directly, we have for every $t\in[0,T]$,
$$\mE\bigg|\int_0^{t}\big[F(r/\eps,X_r^\eps,Y_r^\eps)-  \bar{\bar F}(Y_r^\eps)\big]\cdot(\mI+\nabla_yu)(Y_r^\eps)
\dif r\bigg|^2\leq C_{T}\Big(\eps^{\beta\wedge 1}
+\big[t\cdot\kappa_1(t/\eps)\big]^2\Big).$$
By assumption (\ref{bG}) and applying similar arguments as in (\ref{UI}),  we can obtain that
\begin{align*}
\mE\left(\int_0^t|G(r/\eps,Y_r^\eps)-\bar{G}(Y_r^\eps)|^2_{HS}\dif r\right)\leq C_{T}\Big(\eps^{\b\wedge 1}+\big[t\cdot\hat\kappa_2(t/\eps)\big]^2\Big).
\end{align*}
As a result, we have for every $t\in[0,T]$,
\begin{align*}
\mE\big|Y_{t}^\eps-  \bar Y_{t} \big|^2\leq C_{\lambda,T}\int_0^t\mE|Y_r^\eps-\bar Y_r|^2\dif r+C_{T}\Big(\eps^{\beta\wedge 1}
+\big[t\cdot\kappa_1(t/\eps)\big]^2
+\big[t\cdot\hat\kappa_2(t/\eps)\big]^2\Big).
\end{align*}
By Gronwall's inequality, we obtain the desired result.
\end{proof}

By localization technique, we give:
\begin{proof}[Proof of Corollary \ref{coro}]
For each $n\in \mathbb{N}$, let $\chi_{n}$ be a nonnegative smooth function in $\mathbb{R}^{d_2}$ with $\chi_{n}(y)= 1$ for all $|y|\leq n$ and   $\chi_{n}(y)= 0$ for  $|y|\geq n+1$. Define
$$
 G_{n}(t,y):=G(t, \chi_{n}(y) y), \qquad F_{n}(t,x,y):=F(t,x,y)\chi_{n}(y).
$$	
Then, $G_n$ is bounded and uniformly elliptic. By (\ref{bF0}) and (\ref{bG0}) we have  that for all $T>0$ and $y\in\mR^{d_2}$,
\begin{align*}
\left|\frac{1}{T}\int_{0}^{T}\!\!\int_{\mR^{d_1}}F_n(t,x,y)\mu_t^y(\dif x)\dif t-\bar{\bar F}(y)\chi_{n}(y)\right|\leq C_n\hat\kappa_1(T),
\end{align*}
 and
 \begin{align*}
\frac{1}{T}\int_{0}^{T}\left|G_n(t,y)-{\bar G}(\chi_{n}(y)y)\right|_{HS}^2\dif t\leq C_n \hat\kappa_2^2(T).
\end{align*}
Thus $ F_{n}, G_{n}$ satisfy the conditions in Theorem \ref{main11}. Let $(X^{n,\eps}_{t},Y^{n,\eps}_{t})$ be the solution to SDE \eqref{sde0} with coefficients $ F, G$  replaced by $F_{n}, G_{n}$ respectively.  By Theorem \ref{main11}, we have that for any $T>0$ and $q\geq 1$,
\begin{align*}
\sup_{t\in[0,T]}\mE|Y_t^{n,\eps}-\bar Y^{n}_t|^q\leq C_{n}\Big( \eps^{(\beta\wedge1) q/2}+\sup_{t\in[0,T]}\big[t\cdot\hat\kappa_1(t/\eps)\big]^q+\sup_{t\in[0,T]}
\big[t\cdot\hat\kappa_2(t/\eps)\big]^q\Big),
\end{align*}
where $C_{n}>0$ is independent of $\eps$, and $\bar Y^{n}_t$ satisfies
\begin{align*}
\dif \bar Y^{n}_t=\bar{\bar F}_{n}(\bar Y^{n}_t)\dif t+ \bar G_{n}(\bar Y^{n}_t)\dif W_t^2,\quad\bar Y^{n}_0=y\in\mR^{d_2}
\end{align*}
with
$$
\bar{\bar F}_{n}(y)=\bar{\bar F}(y)\chi_{n}(y),\quad \bar G_{n}(y)={\bar G}(\chi_{n}(y)y).
$$
For each $\eps, n >0$,  define the stopping time  by
$$
\tau_n^{\varepsilon}:=\inf \left\{t \geqslant 0:\left|Y_t^{\varepsilon}\right|+\left|\bar{Y}_t\right| \geqslant n\right\} .
$$
By the uniqueness of the strong solution to SDE \eqref{sde0} and \eqref{bsde} (see e.g. \cite{RX1}), we have
$$
Y_t^{\varepsilon}=Y_t^{n, \varepsilon}\quad\text{and}\quad \bar{Y}_t=\bar{Y}_t^n,\quad \forall \,t \in\left[0, \tau_n^{\varepsilon}\right] .
$$
As a result, we deduce that for any $T>0$ and $0<q<\theta$,
$$
\begin{aligned}
&\sup _{t \in[0, T]}  \mathbb{E}\left|Y_t^{\eps}-\bar{Y}_t\right|^q \leqslant \sup _{t \in[0, T]} \mathbb{E}\left(\left|Y_t^{\eps}-\bar{Y}_t\right|^q \cdot 1_{\left\{t \leqslant \tau_n^{\eps}\right\}}\right)+\sup _{t \in[0, T]} \mathbb{E}\left(\left|Y_t^{\eps}-\bar{Y}_t\right|^q \cdot 1_{\left\{t>\tau_n^{\eps}\right\}}\right) \\
& \leqslant \sup _{t \in[0, T]} \mathbb{E}\left|Y_t^{n, \eps}-\bar{Y}_t^n\right|^q+C\left[\sup _{t \in[0, T]} \mathbb{E}\left(\left|Y_t^{\eps}\right|^\theta+\left|\bar{Y}_t\right|^\theta\right)\right]^{q / \theta}\left[\mathbb{P}\left(T>\tau_n^{\eps}\right)\right]^{(\theta-q) / \theta} \\
& \leqslant {C_{n}} \[\eps^{(\beta\wedge1) q/2}+\sup_{t\in[0,T]}\big[t\cdot\hat\kappa_1(t/\eps)\big]^q+\sup_{t\in[0,T]}
\big[t\cdot\hat\kappa_2(t/\eps)\big]^q\]+C_{\theta,T}  n^{q-\theta},
\end{aligned}
$$
where the last inequality follows by Chebyshev's inequality and  the assumption {\bf($\Theta$)} .  Letting  $\eps \rightarrow 0$ first and then $n \rightarrow \infty$, we  get the desired result.
\end{proof}

\section{Functional central limit theorem}

In this section, we study the normal deviations for the multi-scale SDE (\ref{sde0}). We shall first   derive two weak fluctuation estimates (the functional law of large numbers type and the  functional central limit theorem   type) and a time-oscillation fluctuation estimate    in Subsection 5.1. Then we prove Theorem \ref{main2}
in Subsection 5.2.

\subsection{Weak fluctuation estimates}

Recall that $Z_t^\eps$ is defined by (\ref{zep}). In view of (\ref{sde0}) and (\ref{bsde}), we have
\begin{align*}
Z_t^\eps&=\!\int_0^t\frac{F(r/\eps,X_r^\eps,Y_r^\eps)-\bar{\bar F}(\bar Y_r)}{\eps^{{(1/2)\wedge\vartheta}}}\dif r+\int_0^t\frac{G(Y_r^\eps)- G(\bar Y_r)}{\eps^{{(1/2)\wedge\vartheta}}}\dif W_r^2\\
&=\!\int_0^t\frac{F(r/\eps,X_r^\eps,Y_r^\eps)-\bar{\bar F}( Y_r^\eps)}{\eps^{{(1/2)\wedge\vartheta}}}\dif r+\int_0^t\frac{\bar{\bar F}( Y_r^\eps)-\bar{\bar F}(\bar Y_r)}{\eps^{{(1/2)\wedge\vartheta}}}\dif r+\int_0^t\frac{G(Y_r^\eps)- G(\bar Y_r)}{\eps^{{(1/2)\wedge\vartheta}}}\dif W_r^2.
\end{align*}
To prove the weak
convergence of $Z_t^\eps$, we shall view the process $(X_t^\eps,Y_t^\eps,Z_t^\eps)$ as a whole
stochastic system. Namely, we consider
\begin{equation*}
\left\{ \begin{aligned}
&\dif X^{\eps}_t =\frac{1}{\eps}b(t/\eps, X^{\eps}_t, Y^{\eps}_t)\dif t+\frac{1}{\sqrt{\eps}}\sigma(t/\eps,X_t^\eps,Y_t^\eps)\dif W^{1}_t,\quad\qquad\qquad\qquad\,\; X_0^\eps=x,\\
&\dif Y^{\eps}_t =F(t/\eps,X^{\eps}_t,Y^{\eps}_t)\dif t+G(Y_t^\eps)\dif W^{2}_t,\;\qquad\qquad\qquad\qquad\qquad\qquad\,\,\,\, Y_0^\eps=y,\\
&\dif Z_t^\eps=\frac{F(t/\eps,X_t^\eps,Y_t^\eps)- \bar{\bar F}( Y_t^\eps)}{\eps^{{(1/2)\wedge\vartheta}}}\dif t\no\\
&\quad\qquad+\frac{\bar{\bar F}( Y_t^\eps)-\bar{\bar F}(\bar Y_t)}{\eps^{{(1/2)\wedge\vartheta}}}\dif t+\frac{ G(Y_t^\eps)-  G(\bar Y_t)}{\eps^{{(1/2)\wedge\vartheta}}}\dif W_t^2, \qquad\qquad\qquad\quad \,\,\, Z_0^\eps=0.
\end{aligned} \right.
\end{equation*}
To shorten the notation, we let
\begin{align}\label{333}
\sL_{z}^\eps:=\sL_{z}^\eps(t/\eps,x,y,\bar y)&:=\frac{1}{\eps^{{(1/2)\wedge\vartheta}}}\big[F(t/\eps,x,y)-\bar {\bar F}(y)\big]\cdot\nabla_z+\frac{1}{\eps^{{(1/2)\wedge\vartheta}}}\big[\bar {\bar F}(y)-\bar {\bar F}(\bar y)\big]\cdot\nabla_z\no\\
&\quad+\frac{1}{2}\tr\left(\left(\frac{ G(y)- G(\bar y)}{\eps^{{(1/2)\wedge\vartheta}}}\right)\left(\frac{ G(y)- G(\bar y)}{\eps^{{(1/2)\wedge\vartheta}}}\right)^*\cdot\nabla^2_z\right),
\end{align}
and
\begin{align}\label{34}
\sL_{y,z}^\eps&:=\sL_{y,z}^\eps(y,\bar y):=\tr\left(G(y)\left(\frac{ G(y)- G(\bar y)}{\eps^{{(1/2)\wedge\vartheta}}}\right)^*\cdot\nabla_y\nabla_z\right).
\end{align}
As before, we denote
\begin{align*}
\bar f(t,s,y,z):=\int_{\mR^{d_1}}f(t,s,x,y,z)\mu_t^y(\dif x).
\end{align*}
Note that the function
$$
\hat f(t,s,x,y,z):=f(t,s,x,y,z)-\bar f(t,s,y,z)
$$
always satisfies the centering condition that
\begin{align}\label{cen222}
\int_{\mR^{d_1}}\hat f(t,s,x,y,z)\mu_t^y(\dif x)=0,\quad\forall (t,s,y,z)\in\mR_+^2\times\mR^{2d_2}.
\end{align}
Consider  the following time-inhomogeneous Kolmogorov equation in $\mR_{+}^2\times\mR^{d_1}\times\mR^{2d_2}$:
   \begin{align}\label{psitt}
\partial_t\Psi(t,s,x,y,z)+\sL_0(t,x,y)\Psi(t,s,x,y,z)=-\hat f(t,s,x,y,z),
\end{align}
where $(s,y,z)\in \mR_{+}\times\mR^{2d_2}$ are parameters. We shall denote by
$$\partial_1 \Psi(t,s,x,y,z):=\partial_t \Psi(t,s,x,y,z),\quad\text{and}\quad\partial_2 \Psi(t,s,x,y,z):=\partial_{s} \Psi(t,s,x,y,z).$$
Define
   \begin{align}\label{defps}
 \overline{\delta F\cdot \nabla_z\Psi}(t,s,y,z) :=\int_{\mR^{d_1}}\delta F(t,x,y)\cdot\nabla_z\Psi(t,s,x,y,z)\mu_t^y(\dif x),
  \end{align}
  where
  $$
  \delta F(t,x,y):=F(t,x,y)-\bar F(t,y).
  $$
The following two weak fluctuation estimates  for the integral functional of the process $(X_t^\eps,Y_t^\eps,Z_t^\eps)$ will play an important role in proving Theorem \ref{main2}.

\bl\label{weaf}Let {\bf $(\mathbf{ H_1})$ }and {\bf$(\mathbf{\tilde  H_2})$} hold, $T>0$. Assume that $\sigma, b, F\in C_p^{\alpha/2,\alpha,1+\beta}$ and $G\in C_b^{1+\beta}$ with $\a,\beta\in (0,1]$. Then for every $f\in C_p^{\alpha/2,1+\beta/2,\alpha,1+\beta,2}$ satisfying the centering condition (\ref{cen222})  and $t\in [0,T],$ we have
\begin{align}
\bigg|\mE\bigg(\int_0^t \big[f(r/\eps,r,X_r^{\eps},Y_r^{\eps},Z_r^{\eps})-{\bar f}(r/\eps,r,Y_r^\eps,Z_r^\eps)\big]\dif r\bigg)\bigg|&\leq C_T\,\eps^{\big(1-(\tfrac{1}{2}\wedge\vartheta)\big)\wedge(1+\beta)/2};  \label{we1}
\end{align}
and for $\vartheta\geq 1/2,$
\begin{align}
&\bigg|\mE\bigg(\int_0^t \frac{f(r/\eps,r,X_r^{\eps},Y_r^{\eps},Z_r^{\eps})-{\bar  f}(r/\eps,r,Y_r^\eps,Z_r^\eps)}{\sqrt{\eps}}\dif r\no\\&\qquad\qquad\qquad\qquad\qquad-\int_0^t\overline{\delta F\cdot \nabla_z\Psi}(r/\eps,r,Y_r^\eps,Z_r^\eps) \dif r\bigg)\bigg|\leq C_T\,\eps^{\beta /2};  \label{we3}
\end{align}
where $C_T>0$ is a constant.
\el
\begin{proof} We divide the proof into two steps.

\noindent
{\bf{Step 1.}} We first prove the estimate (\ref{we1}).
 Let $\Psi\in C_p^{1+\alpha/2,1+\b/2,2+\alpha,1+\beta,2}$ be the unique solution of  the time-inhomogeneous Kolmogorov equation (\ref{psitt}), and $\Psi_n$ be the mollifier of $\Psi$ given by
\begin{align}\label{psin}
\Psi_n(t,s,x,y,z):=\Psi*\rho^n:=\int_{\mR^{d_2}}\Psi(t,s,x,y-\tilde{y},z)\rho^{n}(\tilde y)\dif \tilde y.
\end{align}
Applying It\^o's formula to  $\Psi_n(t/\eps,t,X_t^\eps,Y_t^\eps,Z_t^\eps)$, we have
	\begin{align*}
	&\mE[\Psi_n(t/\eps,t,X_t^{\eps},Y_t^\eps,Z_t^\eps)]\\&=\Psi_n(0,0,x,y,0)+
	\frac{1}{\eps}\mE\left(\int_0^t\partial_1\Psi_n(r/\eps,r,X_r^{\eps},Y_r^\eps,Z_r^\eps)\dif r\right)\\
&\quad+\mE\left(\int_0^t(\partial_2+\sL_y)\Psi_n(r/\eps,r,X_r^{\eps},Y_r^\eps,Z_r^\eps)\dif r\right)\\	&\quad+\frac{1}{\eps}\mE\left(\int_0^t\sL_0\Psi_n(r/ \eps,r,X_r^{\eps},Y_r^\eps,Z_r^\eps)\dif r\right)\\
&\quad+\mE\left(\int_0^t(\sL_z^\eps+\sL_{y,z}^\eps)\Psi_n(r/ \eps,r,X_r^{\eps},Y_r^\eps,Z_r^\eps)\dif r\right)
	\end{align*}
where $\sL_0$, $\sL_y$, $\sL_z^\eps$ and  $\sL_{y,z}^\eps$ are defined by (\ref{L0}), (\ref{L1}), (\ref{333}) and (\ref{34}), respectively.
	Multiplying  both sides of the above equality by $\eps$ and taking into account (\ref{psitt}), we obtain
	\begin{align}\label{oi}
\sZ(t,\eps)&:=\bigg|\mE\bigg(\int_0^t \big[f(r/\eps,r,X_r^{\eps},Y_r^{\eps},Z_r^{\eps})-\bar{ f}(r/\eps,r,Y_r^\eps,Z_r^\eps)\big]\dif r\bigg)\bigg|\no\\	&\leq \eps\,\mE\big|\Psi_n(0,0,x,y,0)
-\Psi_n(t/\eps,t,X_t^{\eps},Y_t^\eps,Z_t^\eps)\big|\no\\
&\quad+\eps\mE\bigg|\int_0^t(\partial_2+\sL_y)\Psi_n(r/\eps,r,X_r^{\eps},Y_r^\eps,Z_r^\eps)\dif r\bigg|\no\\	
&\quad+\eps\mE\bigg|\int_0^t(\sL_z^\eps+\sL_{y,z}^\eps)\Psi_n(r/ \eps,r,X_r^{\eps},Y_r^\eps,Z_r^\eps)\dif r\bigg|\no\\
&\quad+\mE\bigg|\int_0^t(\partial_1+\sL_0)(\Psi_n-\Psi)(r/\eps,r,X_r^{\eps},Y_r^\eps,Z_r^\eps)\dif r\bigg|=:\sum_{i=1}^4\sZ_{i}(t,\eps).
	\end{align}
We only need to handle  the term $\sZ_{3}(t,\eps)$, and the other terms can be estimated similarly as in Lemma \ref{xx}.
Since $\Psi\in C_p^{1+\alpha/2,1+\b/2,\alpha+2,1+\beta,2}$,    by definitions (\ref{333}) and (\ref{34}),
we have
\begin{align*}
&\sZ_{3}(t,\eps)
\leq\eps^{1-(\tfrac{1}{2}\wedge\vartheta)}\mE\left|\int_0^t[F(r/\eps,X_r^\eps,Y_r^\eps)-\bar {\bar F}( Y_r^\eps)]\nabla_z\Psi_n(r/ \eps,r,X_r^{\eps},Y_r^\eps,Z_r^\eps)\dif r\right|\no\\
&+\eps\Bigg[\mE\bigg|\int_0^t\!\!\[\frac{\bar {\bar F}(Y_r^\eps)-\bar{\bar F}(\bar Y_r)}{\eps^{(1/2)\wedge\vartheta}}\]\nabla_z\Psi_n(r/ \eps,r,X_r^{\eps},Y_r^\eps,Z_r^\eps)\dif r\bigg|\!\\
&+\frac{1}{2}\mE\!\left|\int_0^t\!\!\tr\left(\left(\frac{ G(Y_r^\eps)- G(\bar Y_r)}{\eps^{(1/2)\wedge\vartheta}}\right)\left(\frac{ G(Y_r^\eps)-G( \bar{Y}_r)}{\eps^{(1/2)\wedge\vartheta}}\right)^*\nabla^2_z\Psi_n(r/ \eps,r,X_r^{\eps},Y_r^\eps,Z_r^\eps)\right)\dif r\right|\!\\
&+\mE\!\left|\int_0^t\!\!\tr\left(G(Y_r^\eps)\left(\frac{ G(Y_r^\eps)- G(\bar Y_r)}{\eps^{(1/2)\wedge\vartheta}}\right)^*
\nabla_y\nabla_z\Psi_n(r/ \eps,r,X_r^{\eps},Y_r^\eps,Z_r^\eps)\right)\dif r\right|\Bigg]\no\\
&=:\sZ_{3,1}(t,\eps)+\sZ_{3,2}(t,\eps).
\end{align*}
It is easy to check that
$$\sZ_{3,1}(t,\eps)\leq C_{T} \eps^{1-(\tfrac{1}{2}\wedge\vartheta)}.$$
For $\sZ_{3,2}(t,\eps),$ due to $b,\sigma\in C_p^{\alpha/2,\alpha,1+\beta},F\in C_p^{\alpha/2,\alpha,1+\beta}$, and by Corollary \ref{bf}, we have $\bar{\bar F}\in C_b^{1+\beta}.$ This together with the condition
$G\in C_b^{1+\beta}$ and the mean value theorem yields that for some $p> 0$ and $C_{T}>0$
\begin{align*}
\sZ_{3,2}(t,\eps)&\leq C_{T} \eps \mE\int_0^t(1+|X_r^\eps|^p)\left[\frac{|\bar {\bar F}(Y_r^\eps)-\bar{\bar F}(\bar Y_r)|}{\eps^{(1/2)\wedge\vartheta}}+\frac{| G(Y_r^\eps)- G(\bar Y_r)|_{HS}^2}{\eps^{1\wedge2\vartheta}}\right]\dif r\no\\
&\quad+C_{T} \eps\, n^{1-(1+\beta)} \mE\int_0^t(1+|X_r^\eps|^{2p})
\left[\frac{| G(Y_r^\eps)- G(\bar Y_r)|_{HS}}{\eps^{(1/2)\wedge\vartheta}}\right]\dif r\no\\
&\leq C_{T} (\eps+\eps n^{1-(1+\beta)})\mE\int_0^t(1+|X_r^\eps|^{4p})(1+|Z_r^\eps|^2)\dif r\no\\&\leq C_{T} (\eps+\eps n^{1-(1+\beta)}).
\end{align*}
Finally, by (\ref{oi}) and taking $n=\eps^{-1/2},$ we obtain
\begin{align*}
\sZ(t,\eps)&\leq C_{T}(\eps+\eps n^{2-(1+\beta)}+\eps^{1-(\tfrac{1}{2}\wedge\vartheta)}+\eps n^{1-(1+\beta)}+n^{-1-\beta})\\
&\leq C_{T}\eps^{\big(1-(\tfrac{1}{2}\wedge\vartheta)\big)\wedge(1+\beta)/2}.
\end{align*}
{\bf{Step 2.}} We proceed to prove the estimate (\ref{we3}).
As in (\ref{oi}) we have
\begin{align}\label{oi11} \sH(t,\eps)&:=\bigg|\mE\bigg(\frac{1}{\sqrt{\eps}}\int_0^t\hat f(r/\eps,r,X_r^{\eps},Y_r^{\eps},Z_r^{\eps})\dif r-\int_0^t\overline{\delta F\cdot\nabla_z\Psi}(r/\eps,r,Y_r^{\eps},Z_r^{\eps})\dif r\bigg)\bigg|\!\no\\	&\leq \sqrt{\eps}\,\mE\big|\Psi_n(0,0,x,y,0)
-\Psi_n(t/\eps,t,X_t^{\eps},Y_t^\eps,Z_t^\eps)\big|\no\\
&\quad+\sqrt{\eps}\mE\bigg|\int_0^t(\partial_2+\sL_y)\Psi_n(r/\eps,r,X_r^{\eps},Y_r^\eps,Z_r^\eps)\dif r\bigg|\no\\	
&\quad+\mE\bigg|\int_0^t(\partial_1+\sL_0+\sqrt{\eps}(\sL_z^\eps+\sL_{y,z}^\eps))(\Psi_n-\Psi)(r/\eps,r,X_r^{\eps},Y_r^\eps,Z_r^\eps)\dif r\bigg|\no\\
&\quad+\bigg|\sqrt{\eps}\mE\int_0^t(\sL_z^\eps+\sL_{y,z}^\eps)\Psi(r/ \eps,r,X_r^{\eps},Y_r^\eps,Z_r^\eps)\dif r\no\\
&\qquad\qquad\qquad\qquad\qquad-\int_0^t\overline{\delta F\cdot\nabla_z\Psi}(r/\eps,r,Y_r^{\eps},Z_r^{\eps})\dif r\bigg|\no\\	
	&\leq C_{T}(\eps^{1/2}+\eps^{1/2} n^{1-\beta}+n^{-1-\beta})\no\\
&
+C_{T}\bigg|\sqrt{\eps}\mE\int_0^t(\sL_z^\eps+\sL_{y,z}^\eps)\Psi(r/ \eps,r,X_r^{\eps},Y_r^\eps,Z_r^\eps)\dif r\no\\
&\qquad\qquad\qquad\qquad\qquad-\int_0^t\overline{\delta F\cdot\nabla_z\Psi}(r/\eps,r,Y_r^{\eps},Z_r^{\eps})\dif r\bigg|,
	\end{align}
where $\Psi_n$ is given by (\ref{psin}).
Now for the last term on the right hand of (\ref{oi11}) , we write
\begin{align*}
&\bigg|\sqrt{\eps}\mE\int_0^t(\sL_z^\eps+\sL_{y,z}^\eps)\Psi(r/ \eps,r,X_r^{\eps},Y_r^\eps,Z_r^\eps)\dif r-\int_0^t\overline{\delta F\cdot\nabla_z\Psi}(r/\eps,r,Y_r^{\eps},Z_r^{\eps})\dif r\bigg|\no\\
&\leq\left|\mE\int_0^t\big[\bar F(r/\eps,Y_r^\eps)-\bar {\bar F}( Y_r^\eps)\big]\nabla_z\Psi(r/ \eps,r,X_r^{\eps},Y_r^\eps,Z_r^\eps)\dif r\right|\no\\
&+\bigg|\mE\int_0^t\big[F(r/\eps,X_r^\eps,Y_r^\eps)-\bar F(r/\eps, Y_r^\eps)\big]\nabla_z\Psi(r/ \eps,r,X_r^{\eps},Y_r^\eps,Z_r^\eps)\dif r\no\\
&\qquad\qquad\qquad\qquad-\int_0^t\overline{\delta F\cdot\nabla_z\Psi}(r/\eps,r,Y_r^{\eps},Z_r^{\eps})\dif r\bigg|\no\\
&+\sqrt{\eps}\bigg[\mE\bigg|\int_0^t\!\!\[\frac{\bar {\bar F}(Y_r^\eps)-\bar{\bar F}(\bar Y_r)}{\sqrt{\eps}}\]\nabla_z\Psi(r/ \eps,r,X_r^{\eps},Y_r^\eps,Z_r^\eps)\dif r\bigg|\!\\
&+\frac{1}{2}\mE\bigg|\int_0^t\!\!\tr\left(\left(\frac{ G(Y_r^\eps)- G(\bar Y_r)}{\sqrt{\eps}}\right)\left(\frac{ G(Y_r^\eps)- G(\bar Y_r)}{\sqrt{\eps}}\right)^*\nabla^2_z\Psi(r/ \eps,r,X_r^{\eps},Y_r^\eps,Z_r^\eps)\right)\dif r\bigg|
\!\\
&+\mE\bigg|\int_0^t\!\!\tr\left(G(Y_r^\eps)\left(\frac{ G(Y_r^\eps)\!-\!G(\bar Y_r)}{\sqrt{\eps}}\right)^*
\!\!\nabla_y\nabla_z\Psi(r/ \eps,r,X_r^{\eps},Y_r^\eps,Z_r^\eps)\right)\dif r\bigg|\bigg]\!\!=:\!\!\sum\limits_{i=1}^3\sH_{i}(t,\eps)
\end{align*}
 Note that $\Psi$ is the solution of the Poisson equation, one can check that the integrands in both $\sH_{1}(t,\eps)$  and $\sH_{2}(t,\eps)$ satisfy the centering condition (\ref{cen222}). Thus, by applying similar arguments as in estimating (\ref{we1}) with $\vartheta\geq 1/2$, we have
$$\sH_{1}(t,\eps)+\sH_{2}(t,\eps)\leq C_{T}\eps^{1/2}.$$
By using similar arguments as in dealing with $\sZ_{3,2}(t,\eps)$ in Step 1, we get
$$\sH_{3}(t,\eps)\leq C_{T}\eps^{1/2}.$$
Substituting these into (\ref{oi11}) and taking $n=\eps^{-1/2},$  we obtain
$$\sH(t,\eps)\leq C_{T}\eps^{\beta /2}. $$
The proof is finished.
\end{proof}

We also provide the following result, which will be used to handle  the time-oscillation terms arising in the proof of Theorem \ref{main2}.
\bl[Time-oscillation fluctuation estimate]\label{time-osc}
Let $T>0,$ and let $h_\varepsilon:\mathbb R_+\times\mathbb R^{d_2}\to\mathbb R^{d_2}$ and $h_0:\mathbb R^{d_2}\to\mathbb R^{d_2}$. Suppose there exists $\eta:\mathbb R_+\to\mathbb R_+$ with $\eta(R)\to0$ as $R\to\infty$ such that for all $t\in[0,T]$ and $y\in\mathbb R^{d_2}$,
\begin{align}\label{time-osc-ass}
\left|
\int_0^t
\big[h_\eps(r/\eps,y)- h_0(y)\big]\dif r
\right|
\leq
t\cdot\eta(t/\eps).
\end{align}
Assume further that for some $\beta\in(0,1]$ and $\varsigma\in[0,1/2]$,
\begin{align}\label{time-osc-reg}
\sup_{t\in[0,T]}\|h_\varepsilon(t,\cdot)\|_{C_b^{1+\beta}} \le C\varepsilon^{-\varsigma},\qquad \|h_0\|_{C_b^{1+\beta}} \le C.
\end{align}
 Let $v\in C_b^{(1+\beta)/2,1+\beta,2}
(\mR_+\times\mR^{d_2}\times\mR^{d_2})$.
Under the assumptions of Theorem \ref{main2}, there exists a constant $C_{T,v}>0$ such that for all $t\in[0,T]$,
\begin{align}\label{time-osc-est}
&\left|
\mE\int_0^t
\big[h_\eps(r/\eps,Y_r^\eps)- h_0(Y_r^\eps)\big]
\cdot v(r,Y_r^\eps,Z_r^\eps)\dif r
\right|
\nonumber\\
&\qquad\leq
C_{T,v}\left(
\eps^{-\varsigma+(1+\beta)/2}
+\sup_{r\in[0,t]}r\cdot\eta(r/\eps)
\right).
\end{align}
The same estimate holds for matrix-valued $h_\eps$ and $ h_0$,
with the Euclidean product replaced by the corresponding trace
product.
\el

\begin{proof}

The proof follows by similar arguments as in dealing with (\ref{Unn}). We only include the details
to clarify the dependence on the factor $\eps^{-\varsigma}$.
Define $H_\varepsilon(t,y):=-\int_0^t[h_\varepsilon(r,y)-h_0(y)]dr$.
Mollifying  $H_\varepsilon$ in the $y$-variable and $v$ in both $t$ and $y$ variables, and using It\^o's
formula gives the integral on the left-hand side of
\eqref{time-osc-est}, together with the boundary term and the
terms involving the generators of $(Y^\eps,Z^\eps)$. The boundary term is controlled by \eqref{time-osc-ass}. The
mollification estimates, the bound \eqref{time-osc-reg}, and the
uniform moment estimates of $(X^\eps,Y^\eps,Z^\eps)$ yield
\begin{align*}
C_{T,v}\left(
\eps^{1-\varsigma} n^{1-\beta}
+n^{-1-\beta}\eps^{-\varsigma}
+\eps^{1-\varsigma}
+\sup_{r\in[0,t]}r\cdot\eta(r/\eps)
\right).
\end{align*}
Taking $n=\eps^{-1/2}$ gives \eqref{time-osc-est}.
The matrix-valued case follows componentwise.
\end{proof}

\subsection{Proof of Theorem \ref{main2}}
  Recall that the limiting processes $\bar Z_t^k,k=1,2,3$ corresponding to Regime 1-3 in Theorem \ref{main2} satisfy (\ref{bze1}), (\ref{bze2}) and (\ref{bze3}), respectively. To prove the weak convergence of $Z_t^\eps$ to $\bar Z_t^k$, we view the  limit processes $\bar Y_t$ and $\bar Z_{t}^k$ as a whole stochastic
system.
Note that the processes $\bar Y_t$ and $\bar Z_{t}^k$ depend on the initial value $y$. Below, we shall write $\bar Y_t(y)$  when we want to stress its dependence on the
initial value, and  use $\bar Z_{t}^k(y,z)$ to denote the process $\bar Z^k_{t}$ with initial point $(\bar Y_0,\bar Z^k_0)=(y,z)\in \mR^{2d_2}.$

Let $\bar \sL_y+\bar \sL^k_{z}+\bar \sL_{y,z}$ be the infinitesimal generator of the  Markov process $(\bar Y_t, \bar Z^k_{t})$,
where  $\bar \sL_y$ is defined by (\ref{bar1}), and
\begin{align}
\bar \sL^1_{z}:=\bar \sL^1_{z}(t,y,z)\!&:=\!\nabla_y\bar {\bar F}(y)z\cdot\!\nabla_z
\!+\!\!\Upsilon_\vartheta(y)\!\cdot\!\nabla_z\!+\!\frac{1}{2}\tr\(\!\!(\nabla_y G(y)z )(\nabla_y G(y)z)^* \! \cdot\!\nabla^2_z\),\label{barz1}
\end{align}
\begin{align}
\bar \sL^2_{z}:=\bar \sL^2_{z}(t,y,z)&:=\nabla_y\bar {\bar F}(y)z\cdot\nabla_z
+\frac{1}{2}
\tr\(\bar{\bar\Sigma}\bar{\bar\Sigma}^*(y)\cdot\nabla^2_z\)\no\\
&\quad\,+\frac{1}{2}\tr\((\nabla_y G(y)z )(\nabla_y G(y)z)^* \cdot\nabla^2_z\),\label{barz3}
\end{align}
\begin{align}
\bar \sL^3_{z}:=\bar \sL^3_{z}(t,y,z)&:=\nabla_y\bar {\bar F}(y)z\cdot\nabla_z+\Upsilon_{1/2}(y)\!\cdot\!\nabla_z
+\frac{1}{2}
\tr\(\bar{\bar\Sigma}\bar{\bar\Sigma}^*(y)\cdot\nabla^2_z\)\no\\
&\quad\,+\frac{1}{2}\tr\((\nabla_y G(y)z )(\nabla_y G(y)z)^* \cdot\nabla^2_z\),\label{barz2}
\end{align}
and
\begin{align}\label{baryz}
\bar \sL_{y,z}:=\bar \sL_{y,z}(y,z):=\tr\( G(y)(\nabla_y G(y)z)^*\cdot\nabla_y\nabla_z\).
\end{align}
Fix $T>0$.  Consider the following backward Kolmogorov equation in $[0,T]\times\mR^{d_2}\times\mR^{d_2}$:
\begin{equation}\left\{\begin{array}{l}\label{PDE}
\displaystyle
\p_t u_k(t,y,z)+(\bar\sL_y+\bar\sL^k_{z}+\bar \sL_{y,z})  u_k(t,y,z)=0,\quad t\in [0, T),\\
u_k(T,y,z)=\varphi(z),
\end{array}\right.
\end{equation}
where $k=1,2,3$, and $\varphi$ is a measurable function on $\mR^{d_2}$. The following probability representation of the solution $u$ of equation (\ref{PDE}) can be proved exactly the same way as in \cite[Theorem 3.4]{RX1}, we omit the details.

\bl
Assume that $\varphi\in C_b^{4}(\mR^{d_2})$. Then  there exists a unique solution $u_k\in C_b^{(2+\beta)/2,2+\beta,4}$ to equation (\ref{PDE}) which is given by
\begin{align}\label{cau}
u_k(t,y,z)=\mE\varphi\big( \bar Z^k_{T-t}(y,z)\big),\quad k=1,2,3.
\end{align}
\el

Now, we are in the position to give:

\begin{proof}[Proof of Theorem \ref{main2}]
	 Given $T>0$ and  $\varphi\in C_b^{4}(\mR^{d_2})$, let $u_k$ with $k=1,2,3$ be defined by (\ref{cau}).
	Then it is obvious that
	$$ u_k(0,y,0)=\mE[\varphi(\bar Z^k_T)].
	$$
	 As a result, we can deduce by It\^o's formula and (\ref{PDE}) that
	\begin{align*}
	&\sR^k(T,\eps):=\left|\mE[\varphi(Z_T^\eps)]-\mE[\varphi( \bar Z^k_T)]\right|=\left|\mE\big[u_k(T,Y_T^\eps,Z_T^\eps)-u_k(0,y,0)\big]\right|\\
	&=\bigg|\mE\left(\int_0^T\big(\p_r+\sL_y+\sL_z^\eps+\sL_{y,z}^\eps\big)u_k(r,Y_r^\eps,Z_r^\eps)\dif r\right)\bigg|\\
	&\leq \bigg|\mE\left(\int_0^T\big(\sL_y-\bar \sL_y\big)u_k(r,Y_r^\eps,Z_r^\eps)\dif r\right)\bigg|\\
	&\quad+\bigg|\mE\left(\int_0^T\big(\sL_{y,z}^\eps-\bar \sL_{y,z}\big)u_k(r,Y_r^\eps,Z_r^\eps)\dif r\right)\bigg|\no\\
&\quad+\bigg|\mE\left(\int_0^T\big(\sL_z^\eps-\bar \sL^k_{z}\big)u_k(r,Y_r^\eps,Z_r^\eps)\dif r\right)\bigg|=:\sum_{i=1}^3\sR^k_i(T,\eps),
	\end{align*}	
		where the operators $\sL_z^\eps,\sL_{y,z}^\eps, \sL_y, \bar \sL_{y,z}$ and $\bar \sL^k_{z}$ are defined by (\ref{333}), (\ref{34}), (\ref{L1}), (\ref{baryz}) and (\ref{barz1})-(\ref{barz2}), respectively.
To control the first term,  by the definitions of $\sL_y$ and $\bar \sL_y$	we have that for $k=1,2,3$,
\begin{align*}
\sR^k_1(T,\eps)&\leq\bigg|\mE\left(\int_0^T\big[F(r/\eps,X_r^\eps,Y_r^\eps)-{\bar F}(r/\eps,Y_r^\eps)\big]\cdot\nabla_y u_k(r,Y_r^\eps,Z_r^\eps)\dif r\right)\bigg|\no\\
&\quad+\bigg|\mE\left(\int_0^T\big[\bar F(r/\eps,Y_r^\eps)-\bar{\bar F}(Y_r^\eps)\big]\cdot\nabla_y u_k(r,Y_r^\eps,Z_r^\eps)\dif r\right)\bigg|\no\\
&=:\sR_{1,1}^k(T,\eps)+\sR_{1,2}^k(T,\eps).
\end{align*}
Note that $[F(t,x,y)-\bar F(t,y)]\cdot\nabla_y u_k(s,y,z)$ satisfies the centering condition (\ref{cen222}). As a direct result of the estimate (\ref{we1}), we have
 $$
 \sR_{1,1}^k(T,\eps)\leq C_{T}\, \eps^{\big(1-(\tfrac{1}{2}\wedge\vartheta)\big)\wedge(1+\beta)/2}.
 $$

 To estimate $\sR^k_{1,2}(T,\eps)$, we apply
Lemma \ref{time-osc} with
$\varsigma=0,
h_\eps(t,y)=\bar F(t,y),
h_0(y)=\bar{\bar F}(y)$ and $
v(s,y,z)=\nabla_yu_k(s,y,z).$
 By assumption (\ref{bF}),  we know that
\begin{align*}
\left|
\int_0^t
\big[\bar F(r/\eps,y)-\bar{\bar F}(y)\big]\dif r
\right|
\leq t\cdot\kappa_1(t/\eps).
\end{align*}
Consequently,
$$\sR^k_{1,2}(T,\eps)\leq C_{T}\,\Big(\eps^{(1+\beta) /2}+\sup\limits_{t\in[0,T]}t\cdot\kappa_1(t/\eps)\Big)\leq C_{T}\,\big(\eps^{(1+\beta) /2}+\eps^{\vartheta}\big),$$
where the last inequality follows from  ${\bf(\tilde H_2)}.$
By the definitions of $\sL^\eps_{y,z}$ and $\bar \sL_{y,z}$, we have for $k=1,2,3$ that
 \begin{align*}
 \sR_2^k(T,\eps)&=\bigg|\mE\bigg(\int_0^T\!\!\tr\bigg(\bigg[ G(Y_r^\eps)\left(\frac{ G(Y_r^\eps)- G(\bar Y_r)}{\eps^{(1/2)\wedge\vartheta}}\right)^*
 \\&\qquad\qquad-G(Y_r^\eps)(\nabla_y G(Y_r^\eps)Z_r^\eps)^*\bigg]\cdot\nabla_y\nabla_z u_k(r,Y_r^\eps,Z_r^\eps)\bigg)\dif r\bigg)\bigg|.
 \end{align*}
Using the fact that $G\in C_b^{1+\beta}$, the mean value theorem and Theorem \ref{main11}, we deduce
 that for some $\theta\in(0,1)$,
	\begin{align}\label{mean} &\sR^k_{2}(T,\eps)\leq\mE\bigg(\int_0^T
\bigg|\tr\bigg(G(Y_r^\eps)\Big(\nabla_yG(Y_r^\eps+\theta(\bar Y_r-Y_r^\eps))\no\\
&\qquad\qquad\quad-\nabla_yG(  Y_r^\eps) \Big) Z_r^\eps\cdot \nabla_y\nabla_z u_k(r, Y_r^\eps,Z_r^{\eps})\bigg)\bigg|\dif r\bigg)\no\\
	&\leq C_{T}\int_0^T\big(\mE|Y_t^\eps-\bar Y_t|^{2\beta}\big)^{1/2}\(\mE\big(1+|Z_t^\eps|\big)^2\)^{1/2}\dif t\leq C_{T}\,\eps^{\beta/2}.
	\end{align}
Below, we proceed to control the last term $\sR^k_{3}(T,\eps)$ in Regime 1-3 separately.

\vspace{1mm}
\noindent
(i) (Regime 1: $\vartheta<1/2$)  In this case, by the definition of the operator $\sL_z^\eps$ in (\ref{333}), we have
\begin{align*}
 &\sR_3^1(T,\eps)\leq
\bigg|\mE\bigg(\int_0^T\!\!\[\frac{\bar {\bar F}(Y_r^\eps)-\bar{\bar F}(\bar Y_r)}{\eps^\vartheta}-\nabla_y\bar {\bar F}( Y_r^\eps)Z_r^\eps\]\cdot\nabla_zu_1(r,Y_r^\eps,Z_r^\eps)\dif r\bigg)\bigg|\\
 &\quad+\frac{1}{2}\bigg|\mE\bigg(\int_0^T\!\!\tr\bigg(\bigg[\left(\frac{ G(Y_r^\eps)- G(\bar Y_r)}{\eps^\vartheta}\right)\left(\frac{ G(Y_r^\eps)- G(\bar Y_r)}{\eps^\vartheta}\right)^*
 \\&\qquad\qquad-(\nabla_y G(Y_r^\eps)Z_r^\eps)(\nabla_y G(Y_r^\eps)Z_r^\eps)^*\bigg]\cdot\nabla^2_z u_1(r,Y_r^\eps,Z_r^\eps)\bigg)\dif r\bigg)\bigg|\\
 &\quad+\bigg|\mE\bigg(\int_0^T\!\!\[\frac{\bar F(r/\eps,Y_r^\eps)-\bar{\bar F}(Y_r^\eps)}{\eps^\vartheta}-\Upsilon_\vartheta(Y_r^\eps)\]\cdot\nabla_z u_1(r,Y_r^\eps,Z_r^\eps)\dif r\bigg)\bigg|\\
 &\quad+\bigg|\mE\bigg(\int_0^T\!\!\[\frac{F(r/\eps,X_r^\eps,Y_r^\eps)-{\bar F}(r/\eps,Y_r^\eps)}{\eps^\vartheta}\]\cdot\nabla_z u_1(r,Y_r^\eps,Z_r^\eps)\dif r\bigg)\bigg|
 =:\!\!\sum_{i=1}^4\sR^1_{3,i}(T,\eps).
 \end{align*}
Using the same arguments as in the estimate of (\ref{mean}), we have
$$\sR^1_{3,1}(T,\eps)+\sR^1_{3,2}(T,\eps)\leq C_{T} \,\eps^{\beta/2}.$$

For $\sR^1_{3,3}(T,\eps)$, set
\begin{align*}
h_\eps(t,y)
&=\frac{\bar F(t,y)-\bar{\bar F}(y)}
{\eps^\vartheta},\quad
h_0(y)=\Upsilon_\vartheta(y),\quad
v(s,y,z)=\nabla_zu_1(s,y,z).
\end{align*}
Assumption ${\bf(\tilde H_3)}$ is precisely the
estimate \eqref{time-osc-ass} with $\varsigma=\vartheta$. Hence,
Lemma \ref{time-osc} gives
\begin{align*}
\sR^1_{3,3}(T,\eps)
\leq C_T\Big(
\eps^{(1+\beta)/2-\vartheta}
+\sup_{t\in[0,T]}t\cdot\tilde\kappa_1(t/\eps)
\Big).
\end{align*}

Note that $[F(t,x,y)-\bar F(t,y)]\nabla_z u_1(s,y,z)$ satisfies the centering condition (\ref{cen222}). Applying (\ref{we1}) directly, we have
$$\sR^1_{3,4}(T,\eps)\leq C_{T}\,\eps^{(1-2\vartheta)\wedge\big((1+\beta)/2-\vartheta\big)}.$$
Combining the above computations, we obtain
$$\sR^1(T,\eps)\leq C_{T}\,\Big(\eps^{\vartheta\wedge(1-2\vartheta)\wedge\beta/2}
+\sup\limits_{t\in[0,T]}t\cdot\tilde\kappa_1(t/\eps)\Big).$$

\noindent
(ii) (Regime 2: $\vartheta>1/2$)
In this case, we have
 \begin{align*}
 \sR_3^2(T,\eps)&\leq
\bigg|\mE\bigg(\int_0^T\!\!\[\frac{\bar {\bar F}(Y_r^\eps)-\bar{\bar F}(\bar Y_r)}{\sqrt{\eps}}-\nabla_y\bar {\bar F}( Y_r^\eps)Z_r^\eps\]\cdot\nabla_zu_2(r,Y_r^\eps,Z_r^\eps)\dif r\bigg)\bigg|\\
 &+\frac{1}{2}\bigg|\mE\bigg(\int_0^T\!\!\tr\bigg(\bigg[\left(\frac{ G(Y_r^\eps)- G(\bar Y_r)}{\sqrt{\eps}}\right)\left(\frac{ G(Y_r^\eps)- G(\bar Y_r)}{\sqrt{\eps}}\right)^*
 \\&\qquad\qquad-(\nabla_y G(Y_r^\eps)Z_r^\eps)(\nabla_y G(Y_r^\eps)Z_r^\eps)^*\bigg]\cdot\nabla^2_z u_2(r,Y_r^\eps,Z_r^\eps)\bigg)\dif r\bigg)\bigg|\\
 &+\bigg|\mE\bigg(\int_0^T\!\!\[\frac{\bar F(r/\eps,Y_r^\eps)-\bar{\bar F}(Y_r^\eps)}{\sqrt{\eps}}\]\cdot\nabla_z u_2(r,Y_r^\eps,Z_r^\eps)\dif r\bigg)\bigg|\\
 &+\bigg|\frac{1}{2}\mE\!\left(\int_0^T\tr\(\bar\Sigma
 \bar\Sigma^*(r/\eps,Y_r^\eps)\cdot\nabla^2_z u_2(r,Y_r^\eps,Z_r^\eps)\)\dif r\right)\\
 &\quad\quad-\frac{1}{2}\mE\!\left(\int_0^T\tr\(\bar{\bar\Sigma}
 \bar{\bar\Sigma}^*(Y_r^\eps)\cdot\nabla^2_z u_2(r,Y_r^\eps,Z_r^\eps)\)\dif r\right)\bigg|\\
  &+\bigg|\mE\bigg(\int_0^T\!\!\[\frac{F(r/\eps,X_r^\eps,Y_r^\eps)-{\bar F}(r/\eps,Y_r^\eps)}{\sqrt{\eps}}\]\cdot\nabla_z u_2(r,Y_r^\eps,Z_r^\eps)\dif r\bigg)\\
 &\quad\quad-\frac{1}{2}\mE\!\left(\int_0^T\tr\(\bar\Sigma
 \bar\Sigma^*(r/\eps,Y_r^\eps)\cdot\nabla^2_z u_2(r,Y_r^\eps,Z_r^\eps)\)\dif r\right)\bigg|=:\!\!\sum_{i=1}^5\sR^2_{3,i}(T,\eps).
 \end{align*}
The first two terms can be estimated exactly the same as before, and we have
$$\sR^2_{3,1}(T,\eps)+\sR^2_{3,2}(T,\eps)\leq C_{T} \,\eps^{\beta/2}.$$
Applying the same arguments as in estimating $\sR^k_{1,2}(T,\eps)$,  by Lemma \ref{time-osc}, assumption (\ref{bF}) and ${\bf(\tilde H_2)}$, we obtain
 $$\sR^2_{3,3}(T,\eps)\leq C_{T}\,\eps^{-1/2}\Big(\eps^{(1+\beta) /2}+\sup\limits_{t\in[0,T]}t\cdot\kappa_1(t/\eps)\Big)\leq  C_{T}\,\eps^{(\vartheta-1/2)\wedge\beta/2}.$$
For $\sR^2_{3,4}(T,\eps)$, we use the matrix-valued version of
Lemma \ref{time-osc} with
$
\varsigma=0,
h_\eps(t,y)=\bar\Sigma\bar\Sigma^*(t,y),
 h_0(y)=\bar{\bar\Sigma}\bar{\bar\Sigma}^*(y)$ and $
v(s,y,z)=\frac12\nabla_z^2u_2(s,y,z).
$
Assumption (\ref{bsigma}) implies that
\begin{align*}
\left|
\int_0^t
\left[
\bar\Sigma\bar\Sigma^*(r/\eps,y)
-\bar{\bar\Sigma}\bar{\bar\Sigma}^*(y)
\right]\dif r
\right|_{HS}
\leq t\cdot\tilde\kappa_2(t/\eps).
\end{align*}
It thus follows that
$$\sR^2_{3,4}(T,\eps)\leq C_{T}\Big(\eps^{(1+\beta) /2}+\sup\limits_{t\in[0,T]}t\cdot\tilde\kappa_2(t/\eps)\Big).$$
To deal with the last term, let
$$
\tilde f(t,s,x,y,z):=F(t,x,y)\cdot\nabla_zu_2(s,y,z)
$$
and
$$
\bar{\tilde f}(t,s,y,z):=\bar F(t,y)\cdot\nabla_zu_2(s,y,z),\quad \hat{\tilde f}(t,s,x,y,z):=\tilde f(t,s,x,y,z)-\bar{\tilde f}(t,s,y,z).
$$
Then one can check that  $\hat{\tilde f}$ satisfies the centering condition (\ref{cen222})
and the solution of the following time-inhomogeneous Kolmogorov equation (where $(s,y,z)$ are parameters)
\begin{equation*}
\partial_t\tilde \psi(t,s,x,y,z)+\sL_0(t,x,y)\tilde \psi(t,s,x,y,z)=-\hat{\tilde f}(t,s,x,y,z)
\end{equation*}
is given by
$$
\tilde \psi(t,s,x,y,z)=\Phi(t,x,y)\cdot\nabla_zu_2(s,y,z),
$$
where $\Phi$ solves the equation (\ref{PPhi}). Thus, by the definitions (\ref{defps}) and (\ref{Sigma}), we have
\begin{align*}
\overline{\delta F\cdot \nabla_z\tilde\psi}(t,s,y,z) &=\int_{\mR^{d_1}}\delta F(t,x,y)\cdot\nabla_z\tilde\psi(t,s,x,y,z)\mu_t^y(\dif x)\\
&=\frac{1}{2}\tr\(\bar\Sigma
\bar\Sigma^*(t,y)\cdot\nabla^2_z u_2(s,y,z)\).
\end{align*}
Consequently, as  a direct result of (\ref{we3}) of Lemma \ref{weaf}, we get
$$\sR^2_{3,5}(t,\eps)\leq C_{T} \eps^{\beta/2} .$$
Combining the above computations, we obtain
$$\sR^2(T,\eps)\leq C_{T}\Big(\eps^{(\vartheta-1/2)\wedge\beta /2}+\sup\limits_{t\in[0,T]}t\cdot\tilde\kappa_2(t/\eps)\Big).$$

\noindent
(iii) (Regime 3: $\vartheta=1/2$) In this case, we have
\begin{align*}
 \sR_3^3(T,\eps)&\leq
\bigg|\mE\bigg(\int_0^T\!\!\[\frac{\bar {\bar F}(Y_r^\eps)-\bar{\bar F}(\bar Y_r)}{\sqrt{\eps}}-\nabla_y\bar {\bar F}( Y_r^\eps)Z_r^\eps\]\cdot\nabla_zu_3(r,Y_r^\eps,Z_r^\eps)\dif r\bigg)\bigg|\\
 &+\frac{1}{2}\bigg|\mE\bigg(\int_0^T\!\!\tr\bigg(\bigg[\left(\frac{ G(Y_r^\eps)- G(\bar Y_r)}{\sqrt{\eps}}\right)\left(\frac{ G(Y_r^\eps)- G(\bar Y_r)}{\sqrt{\eps}}\right)^*
 \\&\qquad\qquad-(\nabla_y G(Y_r^\eps)Z_r^\eps)(\nabla_y G(Y_r^\eps)Z_r^\eps)^*\bigg]\cdot\nabla^2_z u_3(r,Y_r^\eps,Z_r^\eps)\bigg)\dif r\bigg)\bigg|\\
 &+\bigg|\mE\bigg(\int_0^T\!\!\[\frac{\bar F(r/\eps,Y_r^\eps)-\bar{\bar F}(Y_r^\eps)}{\sqrt{\eps}}-\Upsilon_{1/2}(Y_r^\eps)\]\cdot\nabla_z u_3(r,Y_r^\eps,Z_r^\eps)\dif r\bigg)\bigg|\\
  &+\bigg|\frac{1}{2}\mE\!\left(\int_0^T\tr\(\bar\Sigma
 \bar\Sigma^*(r/\eps,Y_r^\eps)\cdot\nabla^2_z u_3(r,Y_r^\eps,Z_r^\eps)\)\dif r\right)\\
 &\quad\quad-\frac{1}{2}\mE\!\left(\int_0^T\tr\(\bar{\bar\Sigma}
 \bar{\bar\Sigma}^*(Y_r^\eps)\cdot\nabla^2_z u_3(r,Y_r^\eps,Z_r^\eps)\)\dif r\right)\bigg|\\
 &+\bigg|\mE\bigg(\int_0^T\!\!\[\frac{F(r/\eps,X_r^\eps,Y_r^\eps)-{\bar F}(r/\eps,Y_r^\eps)}{\sqrt{\eps}}\]\cdot\nabla_z u_3(r,Y_r^\eps,Z_r^\eps)\dif r\bigg)\\
 &\quad\quad-\frac{1}{2}\mE\!\left(\int_0^T\tr\(\bar\Sigma
 \bar\Sigma^*(r/\eps,Y_r^\eps)\cdot\nabla^2_z u_3(r,Y_r^\eps,Z_r^\eps)\)\dif r\right)\bigg|=:\!\!\sum_{i=1}^5\sR^3_{3,i}(T,\eps).
 \end{align*}
The terms $\sR^3_{3,1}(T,\eps)$,  $\sR^3_{3,2}(T,\eps)$, $\sR^3_{3,4}(T,\eps)$ and $\sR^3_{3,5}(T,\eps)$ can be controlled the same as in Regime 2, and we have
$$\sR^3_{3,1}(T,\eps)+\sR^3_{3,2}(T,\eps)+\sR^3_{3,5}(T,\eps)\leq C_{T} \,\eps^{\beta/2},$$
and
$$\sR^3_{3,4}(T,\eps)\leq C_{T}\Big(\eps^{(1+\beta) /2}+\sup\limits_{t\in[0,T]}t\cdot\tilde\kappa_2(t/\eps)\Big).$$

For $\sR^3_{3,3}$, we take
$\varsigma=1/2$,
\begin{align*}
h_\eps(t,y)
=\frac{\bar F(t,y)-\bar{\bar F}(y)}{\sqrt\eps},
\qquad
 h_0(y)=\Upsilon_{1/2}(y),
\end{align*}
in Lemma \ref{time-osc} and use ${\bf(\tilde H_3)}$ to deduce that
\begin{align*}
\sR^3_{3,3}(T,\eps)
\leq C_T\Big(
\eps^{\beta/2}
+\sup_{t\in[0,T]}t\cdot\tilde\kappa_1(t/\eps)
\Big).
\end{align*}
Combining the above computations, we obtain
$$\sR^3(T,\eps)\leq C_{T}\,\Big(\eps^{\beta/2}
+\sup\limits_{t\in[0,T]}t\cdot\tilde\kappa_1(t/\eps)
+\sup\limits_{t\in[0,T]}t\cdot\tilde\kappa_2(t/\eps)\Big).$$
The proof is finished.
\end{proof}




\end{document}